\documentclass[a4paper]{amsart}
\usepackage[all]{xy}
\usepackage{amssymb}
\usepackage{xypic}
\usepackage{array}
\usepackage{tikz}
\usepackage{tabularx}
\usepackage{color}
\usepackage{hyperref}
\usepackage{float}
\numberwithin{equation}{section}

\newtheorem{theorem}{Theorem}[section]
\newtheorem{lemma}[theorem]{Lemma}
\newtheorem{proposition}[theorem]{Proposition}
\newtheorem{corollary}[theorem]{Corollary}
\newtheorem{definition}[theorem]{Definition}

\theoremstyle{remark}
\newtheorem{remark}[theorem]{Remark}

\newcommand{\bpr}{\begin{proposition}}
\newcommand{\epr}{\end{proposition}}
\newcommand{\btm}{\begin{theorem}}
\newcommand{\etm}{\end{theorem}}
\newcommand{\bco}{\begin{corollary}}
\newcommand{\eco}{\end{corollary}}
\newcommand{\blm}{\begin{lemma}}
\newcommand{\elm}{\end{lemma}}
\newcommand{\bdf}{\begin{definition}}
\newcommand{\edf}{\end{definition}}
\newcommand{\bpm}{\begin{pmatrix}}
\newcommand{\epm}{\end{pmatrix}}
\newcommand{\beq}{\begin{equation}}
\newcommand{\eeq}{\end{equation}}
\newcommand{\bit}{\begin{itemize}}
\newcommand{\eit}{\end{itemize}}
\newcommand{\brm}{\begin{remark}}
\newcommand{\erm}{\end{remark}}
\newcommand{\bpf}{\begin{proof}}
\newcommand{\epf}{\end{proof}}

\newcommand{\lie}{\mathfrak}
\newcommand{\nr}{\textnormal}
\newcommand{\wt}{\widetilde}

\newcommand{\ol}{\overline}

\renewcommand{\AA}{\mathbb{A}}

\newcommand{\CC}{\mathbb{C}}

\newcommand{\MM}{\mathbb{M}}
\newcommand{\NN}{\mathbb{N}}
\newcommand{\PP}{\mathbb{P}}

\newcommand{\RR}{\mathbb{R}}

\newcommand{\Ee}{\mathcal{E}}

\newcommand{\Jj}{\mathcal{J}}

\newcommand{\Mm}{\mathcal{M}}
\newcommand{\Nn}{\mathcal{N}}
\newcommand{\Oo}{\mathcal{O}}

\newcommand{\Qq}{\mathcal{Q}}

\newcommand{\gG}{\lie{g}}

\newcommand{\zZ}{\lie{z}}

\newcommand{\A}{\mathbf{A}}
\newcommand{\B}{\mathbf{B}}

\newcommand{\E}{\mathbf{E}}
\newcommand{\F}{\mathbf{F}}

\newcommand{\I}{\mathbf{I}}
\newcommand{\J}{\mathbf{J}}

\newcommand{\V}{\mathbf{V}}
\newcommand{\W}{\mathbf{W}}
\newcommand{\X}{\mathbf{X}}

\newcommand{\Tr}{\mathbf{tr}}

\newcommand{\g}{\mathbf{g}}
\newcommand{\h}{\mathbf{h}}
\renewcommand{\k}{\mathbf{k}}

\newcommand{\quotient}[2]{{\raisebox{.2em}{\thinspace $#1$}\left / \raisebox{-.15em}{ $#2$}\right.}}
\newcommand{\git}[2]{{\raisebox{.2em}{\thinspace $#1$}\left /\!\!/ \raisebox{-.15em}{$#2$}\right.}}
\newcommand{\morph}[6]{\begin{array}{cccc} #1 \quad : &  #2  & \stackrel{#6}{\longrightarrow} &  #3  \\  & #4 &\longmapsto & #5  \end{array}}

\newcommand{\gitchar}[3]{{\raisebox{.2em}{\thinspace $#1$}\left / \!\! /_{#3} \raisebox{-.15em}{\thinspace $#2$}\right.}}

\newcommand{\qua}{\thinspace}

\newcommand{\lra}{\longrightarrow}

\newcommand{\id}{\mathbf{1}}

\DeclareMathOperator{\im}{im}

\DeclareMathOperator{\tr}{tr}

\DeclareMathOperator{\head}{head}
\DeclareMathOperator{\tail}{tail}

\DeclareMathOperator{\GL}{GL}

\DeclareMathOperator{\SL}{SL}
\DeclareMathOperator{\PSL}{PSL}
\DeclareMathOperator{\U}{U}
\DeclareMathOperator{\SU}{SU}
\DeclareMathOperator{\PSU}{PSU}
\DeclareMathOperator{\Ort}{O}

\DeclareMathOperator{\Sp}{Sp}

\DeclareMathOperator{\Hom}{Hom}

\DeclareMathOperator{\Aut}{Aut}

\title{Branes in the moduli space of framed sheaves}

\author{Emilio Franco}
\address{Emilio Franco \\ IMECC (Instituto de Matem\'atica, Estat\'istica e Computa\c{c}\~ao Cien\-t\'i\-fi\-ca) \\ Universidade Estadual de Campinas \\ 
Rua S\'ergio Buarque de Holanda 651 \\ Cidade Universit\'aria "Zeferino Vaz", Campinas (SP, Brazil)}
\email{emilio\_franco@ime.unicamp.br}

\author{Marcos Jardim}
\address{Marcos Jardim \\ IMECC (Instituto de Matem\'atica, Estat\'istica e Computa\c{c}\~ao Cien\-t\'i\-fi\-ca) \\ Universidade Estadual de Campinas \\ 
Rua S\'ergio Buarque de Holanda 651 \\ Cidade Universit\'aria "Zeferino Vaz", Campinas (SP, Brazil)}
\email{jardim@ime.unicamp.br}

\author{Simone Marchesi}
\address{Simone Marchesi \\ IMECC (Instituto de Matem\'atica, Estat\'istica e Computa\c{c}\~ao Cien\-t\'i\-fi\-ca) \\ Universidade Estadual de Campinas \\ 
Rua S\'ergio Buarque de Holanda 651 \\ Cidade Universit\'aria "Zeferino Vaz", Campinas (SP, Brazil)}
\email{marchesi@ime.unicamp.br}

\begin{document}

\thanks{EF is supported by the FAPESP postdoctoral grant number 2012/16356-6 and and FAPESP BEPE grant number 2015/06696-2. MJ is partially supported by the CNPq grant number 303332/2014-0 and the FAPESP grant number 2014/14743-8. SM is partially supported by the FAPESP grant number 2014/19676-7.}

\date{\today}

\maketitle

\begin{abstract}
In the physicist's language, a {\it brane} in a hyperk\"ahler manifold is a submanifold which is either complex or lagrangian with respect to three K\"ahler structures of the ambient manifold. By considering the fixed loci of certain involutions, we describe branes in Nakajima quiver varieties of all possible types. We then focus on the moduli space of framed torsion free sheaves on the projective plane, showing how the involutions considered act on sheaves, and proving the existence of branes in some cases.
\end{abstract}

\tableofcontents


\section{Introduction}

Given a K\"ahler manifold $Y$ with complex structure $\Gamma$, Riemannian metric $\eta$ and symplectic form  
\beq
\label{eq definition of omega_i}
\omega(\qua \cdot \qua, \qua \cdot \qua) := \eta(\qua \cdot \qua, \Gamma(\qua \cdot \qua)),
\eeq
we say that submanifold $Y'\subset Y$ is a $A$-{\it brane} if it is lagrangian with respect to $\omega$, and that it is a $B$-{\it brane} if it is a complex submanifold.

Now let $Y$ be a hyperk\"ahler manifold with complex structures $\Gamma_1$, $\Gamma_2$, and $\Gamma_3$ satisfying the usual quaternionic relations, and Riemannian metric $\eta$; denote by $\omega_1$, $\omega_2$ and $\omega_3$ the associated symplectic forms, constructed as in \eqref{eq definition of omega_i}. A subvariety $Y'\subset Y$ is said to be a {\it brane} if it is either an $A$-brane (i.e. lagrangian) or a $B$-brane (complex) with respect to each symplectic form or complex structure. One then specifies the behaviour of $Y'$ by saying that $Y'$ is either a $(B,B,B)$, $(B,A,A)$, $(A,B,A)$ or $(A,A,B)$-brane; note that those are all the possible branes. A $(B,B,B)$-brane is a hyperk\"ahler submanifold of $Y$, since the hyperk\"ahler structure of $Y$ restricts to the brane. On the other hand, a $(B,A,A)$-brane is a complex lagrangian submanifold with respect to the holomorphic symplectic form $\Omega_1 = \omega_2 + i \omega_3$. Similarly, $(A,B,A)$ and $(A,A,B)$-branes correspond with complex lagrangian 
submanifolds with respect to $\Omega_2$ and $\Omega_3$. 

In the groundbreaking paper of Kapustin and Witten \cite{KW}, branes in the moduli space of Higgs bundles played an important role in connection with the geometric Langlands program and mirror symmetry. This explains the great interest that the topic has attracted in the recent years, achieving the description of branes within moduli space of Higgs bundles as fixed point loci of involutions associated with reductions of structure group, and real structures \cite{baraglia&schaposnik_1, baraglia&schaposnik_2, biswas&oscar&hurtubise, oscar&biswas, oscar_1, oscar_2, hitchin}. 

In this paper, we study branes in Nakajima quiver varieties, another important class of hyperk\"ahler manifolds introduced in \cite{nakajima_1}. Differently from the case of Higgs bundles, our motivation is purely geometrical, focusing on the study of various involutions on these varieties and, in the special case of the Jordan quiver, their connection with interesting classes of locally free sheaves on $\PP^2$.

Indeed, Nakajima quiver varieties are constructed as hyperk\"ahler quotients of the vector space of representations of a given quiver. Using linear algebra transformations (transposition, multiplication by scalars, addition and conjugation) we construct involutions on this space representations. Since these involutions are compatible with the hyperk\"ahler structure (they preserve the Riemannian metric and either commute or anticommute with the three complex structures), and with the action of the group (they preserve the orbits and the preimages of the moment maps), one obtains involutions in the quotient that are compatible with the hyperk\"ahler structure. If nonempty, the fixed point loci of these involutions provide examples of all possible kinds of branes on Nakajima quiver varieties. 

We then focus on the study of the Nakajima quiver variety for the Jordan quiver, which is coincides, via de ADHM correspondence, with the moduli space of framed torsion free sheaves over $\PP^2$. We describe the involutions previously obtained in terms of transformations of the framed torsion free sheaves:  

\begin{itemize}
\item In the first case, the involution is described by means of the dual framed sheaf, so one has to restrict to the locally free locus of the moduli space. The fixed points are framed autodual bundles, symplectic and orthogonal, which have been widely studied in the literature, see for example \cite{AB, CO,  marcos&simone&anna, O, scalise}. We show that the autodual locus is a $(B,B,B)$-brane inside the moduli space of framed torsion free sheaves over $\PP^2$. 

\item The following involutions considered can be understood in terms of pull-backs by unitary holomorphic involutions in $\PP^2$ that preserve the line at infinite. Choosing projective coordinates $[x_0 : x_1 : x_2]$ such that the line at infinity is $x_0 = 0$, a unitary involution is constructed with an involutive element $M$ of the unitary group $\U(2)$ acting on $x_1$ and $x_2$. The fixed point set is given by those framed locally free sheaves which are isomorphic to the pull-back by the unitary involution. If $\det M = 1$, this gives a $(B,B,B)$-brane, while if $\det M = -1$, one has a $(B,A,A)$-brane. 

\item The last involutions considered are given by the composition of the complex conjugation of the framed sheaf and the pull-back by an antiholomorphic (unitary) involution of $\PP^2$. Note that the composition of these two operations is a functor within the holomorphic category. In this case, if $\det M = 1$, the fixed point set is a $(A,B,A)$-brane, and if $\det M = -1$, it is a $(A,A,B)$-brane.
\end{itemize}

We also prove the existence of fixed points for some particular values of the rank and the charge, thus showing the nonemptiness of the branes considered.


\section{Preliminaries}

We begin by revising the construction of Nakajima quiver varieties and describing how to construct branes on hyperk\"ahler quotients using involution on the covering hyperk\"ahler manifold.


\subsection{The hyperk\"ahler structure of Nakajima quiver varieties}
\label{sc hyperkahler structure of Nn}

In this section we provide a review of the construction of Nakajima quiver varieties. We address the reader to Nakajima's original work \cite{nakajima_1} for further details.

Let $\Qq$ be a quiver with set of vertices $\Delta$ and set of arrows $\Sigma$. Define the {\it framed quiver} $\Qq^\diamond$ as the quiver whose vertex set is $\Delta \sqcup \Delta'$, where $\Delta'$ is another copy of the vertex set $\Delta$, and whose arrow set is $\Sigma^\diamond = \Sigma \sqcup \Sigma_\Delta$, where $\Sigma_\Delta$ is a set of arrows starting on the vertex $i$ of $\Delta$ and ending on the corresponding vertex $i' \in \Delta'$, for each $i \in \Delta$. 

For each arrow $a \in \Sigma$ with head $\head(a) \in \Delta$ and tail $\tail(a) \in \Delta$, construct another arrow $\wt{a}$ with head $\tail(a)$ and tail $\head(a)$. Consider the {\it opposite quiver} $\Qq^\vee$ to be the quiver with the same set of vertices $\Delta$ and whose set of arrows is $\Sigma^\vee = \{ \wt{a} : a \in \Sigma \}$. We also define the {\it double quiver} $\wt{\Qq}$ as the quiver with vertex set $\Delta$ and arrow set $\Sigma \sqcup \Sigma^\vee$. 

Let $\V$ and $\W$ be two collections of hermitian vector spaces defined as $\V = \{ V_i \}_{i \in \Delta}$ and $\W = \{ W_i \}_{i \in \Delta}$. Let us define the following vector spaces
\[
\Hom^\Sigma(\V, \W) = \bigoplus_{a \in \Sigma} \Hom(V_{\tail(a)}, W_{\head(a)})
\]
and 
\[
\Hom^\Delta(\V, \W) = \bigoplus_{i \in \Delta} \Hom(V_i, W_i).
\]
We call $\dim(\V) = (\dim(V_i))_{i \in \Delta} \in \NN^{| \Delta |}$ the \emph{dimension vector} of the collection $\V$.

A \emph{representation} for the quiver $\Qq$ associated to the collection $\V$ is an element of $\Hom^\Sigma(\V,\V)$. Given a representation $\A \in \Hom^\Sigma(\V, \V)$, set  
\[
\A^* := ( A_a^* )_{\wt{a} \in \Sigma^\vee} \in \Hom^{\Sigma^\vee}(\V^*, \V^*),
\]
\[
\A^\vee := ( A_a^\vee )_{\wt{a} \in \Sigma^\vee} \in \Hom^{\Sigma^\vee}(\V^\vee, \V^\vee)
\]
and 
\[
\ol{\A} := ( \ol{A}_a )_{a \in \Sigma} \in \Hom^{\Sigma}(\ol{\V},\ol{\V}),
\]
where $(\cdot)^*$ denotes the hermitian adjoint operator, $(\cdot)^\vee$ the transposition operator and $\overline{(\cdot)}$ the complex conjugate.

Given $\A \in \Hom^{\Sigma}(\V, \V')$, $\B \in \Hom^{\Sigma^\vee}(\V', \V'')$, $\J \in \Hom^{\Delta}(\V, \V')$ and $\I \in \Hom^{\Delta}(\V', \V'')$, define 
\[
\B \A = \left ( \sum_{\tail(a) = i} B_{\wt{a}}A_a \right )_{i \in \Delta} \in \Hom^{\Delta}(\V, \V''), 
\]
and
\[
\I \J = \left ( I_i J_i \right )_{i \in \Delta} \in \Hom^{\Delta}(\V, \V'').
\]

Note that the quiver $\wt{\Qq^\diamond}$ has as vertex set $\Delta \sqcup \Delta'$ and arrow set $\wt{\Sigma^\diamond} = \Sigma \sqcup \Sigma^\vee \sqcup \Sigma_\Delta \sqcup \Sigma_\Delta^\vee$. 
Denote the space of representations of $\wt{\Qq^\diamond}$ for the collections of vectors spaces $(\V, \W)$ by
\[
\MM^\Qq_{\V, \W} := \Hom^{\Sigma}(\V, \V) \oplus \Hom^{\Sigma^\vee}(\V, \V) \oplus \Hom^{\Delta}(\W, \V) \oplus \Hom^{\Delta}(\V, \W).
\]
Thus a representation $\X$ of $\wt{Q^\diamond}$ for the collections of vector spaces $(\V, \W)$, is a quadruple $(\A, \B, \I, \J)$, where
\[
\A \in \Hom^{\Sigma}(\V, \V) ~~,~~ \B \in \Hom^{\Sigma^\vee}(\V, \V), 
\]
\[
\I \in \Hom^{\Delta}(\W, \V) ~~{\rm and}~~ \J \in \Hom^{\Delta}(\V, \W).
\]
 
Define the \emph{trace} of $\A = ( A_i )_{i \in \Delta} \in \Hom^\Delta(\V, \V)$ as
\[
\Tr \A := \sum_{i \in \Delta} \tr A_i. 
\]
 
One can endow $\MM^\Qq_{\V,\W}$ with a natural (flat) hyperk\"ahler structure in the following way. First, we consider the hermitian metric $\eta$ on $\MM^\Qq_{\V,\W}$ given by
\begin{align*}
\eta(\X,\X') := &\frac{1}{2} \Tr\left( \A (\A')^* + \A' \A^* \right) +
\frac{1}{2} \Tr\left( \B (\B')^* + \B' \B^* \right) + 
\\
&+ \frac{1}{2} \Tr\left( \I (\I')^* + \I' \I^* \right) + 
\frac{1}{2} \Tr\left( \J (\J')^* + \J' \J^* \right),
\end{align*}
where $\X = (\A, \B, \I, \J)$ and $\X' = (\A', \B', \I', \J')$. Next, consider the following complex structures on $\MM^\Qq_{\V,\W}$
\beq
\label{eq complex structures}
\begin{array}{lll}
         \Gamma_1(\A, \B, \I, \J) := (i\A, i\B, i\I, i\J),          \\

         \Gamma_2(\A, \B, \I, \J) := (-\B^*, \A^*, -\J^*, \I^*),        \\

         \Gamma_3(\A, \B, \I, \J) := (-i\B^*, i\A^*, -i\J^*, i\I^*).
        \end{array}
\eeq
Observe that the $\Gamma_k$ satisfy the quaternionic relations. Also, note that our metric is compatible with these complex structures,
\[
\eta(\Gamma_k(\X), \Gamma_k(\X')) = \eta(\X,\X'),
\]
and therefore, using \eqref{eq definition of omega_i}, one can define the symplectic forms on $\omega_1$, $\omega_2$ and $\omega_3$, which completes the description of $\MM^\Qq_{\V,\W}$ as a hyperk\"ahler manifold. 

Next, consider the groups, defined by the unitary and general linear groups,
\[
\U(\V) := \Pi_{i \in \Delta} \U(V_i)
\]
and 
\[
\GL(\V) := \Pi_{i \in \Delta} \GL(V_i).
\]
Both groups act on $\MM^\Qq_{\V,\W}$ in the following manner; given an element $\g = ( g_i )_{i \in \Delta}$ of $\GL(\V)$ or $\U(\V)$, define:
\[
\A\g := (A_i g_{\tail(a)})_{a \in \Sigma},
\]
\[
\g \B := (g_{\head(a)} B_i)_{a \in \Sigma},
\]
\[
\g \I := (g_i I_i)_{i \in \Delta},
\] 
and
\[
\J \g := (J_i g_i)_{i \in \Delta}.
\]
One can take the following action of $\GL(\V)$ or $\U(\V)$ on $\MM^\Qq_{\V,\W}$ by setting
\beq
\label{eq definition of the action of U V on MM}
\g \cdot (\A, \B, \I, \J) := ( \g\A\g^{-1}, \g\B\g^{-1}, \g\I, \J\g^{-1}).
\eeq
One easily checks that for any $\X \in \MM^\Qq_{\V,\W}$ we have
\[
\Gamma_k(\g \cdot \X) = \g \cdot \Gamma_k(\X) ~~ {\rm and} ~~ 
\eta(\g \cdot \X, \g \cdot \X) = \eta(\X,\X'), 
\]
so this action preserves the three complex structures and the hermitian metric. It follows that the action preserves the three symplectic forms,
\[
\omega_k(\g \cdot \X, \g \cdot \X) = \omega_k(\X,\X'). 
\]

Analogously, we have the action of $\h \in \U(\W)$ on $\MM^\Qq_{\V,\W}$ given by 
\beq
\label{eq definition of the action of U W on MM}
\h \cdot (\A, \B, \I, \J) := ( \A, \B, \I \h, \h^{-1}\J).
\eeq
As before, this action commutes with the three complex structures $\Gamma_k$, preserves the metric $\eta$, and therefore, the three symplectic forms $\omega_k$. One can easily see that both actions \eqref{eq definition of the action of U V on MM} and \eqref{eq definition of the action of U W on MM} commute, so it is possible to define an action of $\U(\V) \times \U(\W)$ on $\MM^\Qq_{\V,\W}$ by setting 
\beq
\label{eq definition of de action of U V U W on MM}
(\g,\h) \cdot \X := ( \g \A \g^{-1}, \g \B \g^{-1}, \g\I \h, \h^{-1}\J\g^{-1}).
\eeq

The action of $\U(\V)$ given by \eqref{eq definition of the action of U V on MM} on the symplectic manifolds $(\MM^\Qq_{\V,\W}, \omega_1)$, $(\MM^\Qq_{\V,\W}, \omega_2)$ and $(\MM^\Qq_{\V,\W}, \omega_3)$ defines the following moment maps (see \cite{nakajima_1} for instance)
\[
\mu_1(\X) = \frac{-1}{2} \left( [\A, \B] + [\A^*, \B^*] + \I \J - \J^* \I^* \right),
\]
\[
\mu_2(\X) = \frac{-1}{2i} \left( [\A, \B] - [\A^*, \B^*] + \I \J + \J^* \I^* \right)
\]
and 
\[
\mu_3(\X) = \frac{i}{2} \left( [\A, \A^*] + [\B, \B^*] + \I \I^* - \J^* \J \right ),
\]
where 
\[
[\A, \B] = \A \B - \B \A = \left (  \sum_{\tail(a) = i} A_{\wt{a}}B_a - B_{\wt{a}} A_a  \right )_{i \in \Delta}. 
\]

Note that $\mu_1$ and $\mu_2$ can be recombined into
\[
\mu_{\CC}(\X) = -\mu_1(\X) -i \mu_2(\X) = [\A, \B] + \I \J.
\]

One checks that the action of $\U(\V) \times \U(\W)$ preserves $\mu_k^{-1}(0)$ for  $k = 1,2,3$; in fact
\[
\mu_k((\g,\h) \cdot \X) = \g \mu_k(\X) \g^{-1}.
\]

Then, $\GL(\V)$ preserves $\mu_\CC^{-1}(0)$ and so, the {\it affine Nakajima quiver variety} for $\Qq$ is defined as the affine GIT quotient  
\beq
\label{eq definition of Nn_0 as a affine GIT quotient}
\Nn^\Qq_0 :=  \git{\mu_\CC^{-1}(0)}{\GL(\V)}.
\eeq

Take the character 
\beq \label{eq definition of chi}
\morph{\chi}{\GL(\V)}{\CC^*}{\g = (g_i)_{i \in \Delta}}{\Pi_{i \in \Delta} \det(g_i).}{} 
\eeq
The {\it Nakajima quiver variety for the quiver $\Qq$} is defined as the GIT quotient
\beq
\label{eq definition of Nn_1 as a projective GIT quotient}
\Nn^\Qq_1 :=  \gitchar{\mu_\CC^{-1}(0)}{\GL(\V)}{\chi}.
\eeq
Using the inverse character, we define analogously
\beq
\label{eq definition of Nn_-1 as a projective GIT quotient}
\Nn^\Qq_{-1} :=  \gitchar{\mu_\CC^{-1}(0)}{\GL(\V)}{\chi^{-1}}.
\eeq
By construction, there are projective morphisms $\Nn^\Qq_1 \to \Nn^\Qq_0$ and $\Nn^\Qq_{-1} \to \Nn^\Qq_0$.

Following \cite{king}, we say that $\X = (\A, \B, \I, \J) \in \mu_\CC^{-1}(0)$ is {\it stable} if there is no proper collection of subspaces $\V' \subset \V$ such that $\A (\V') \subset \V'$, $\B(\V') \subset \V'$ and $\im \I \subset \V'$. Similarly, $\X$ is {\it costable} if there is no proper collection of subspaces $\V' \subset \V$ such that $\A (\V') \subset \V'$, $\B(\V') \subset \V'$ and $\V' \subset \ker \J$. We say that $\X$ is {\it regular} if it is both stable and costable. Denote by $\mu_\CC^{-1}(0)^{st}$, $\mu_\CC^{-1}(0)^{ct}$ and $\mu_\CC^{-1}(0)^{reg}$ the $\GL(\V)$-invariant sets of stable, costable and regular points.

Semistability, polystability and stability coincide in the GIT quotients \eqref{eq definition of Nn_1 as a projective GIT quotient} and \eqref{eq definition of Nn_-1 as a projective GIT quotient}. In this case, \cite[Proposition 3.1]{king} implies the following.

\bpr 
A point $\X \in \mu_\CC^{-1}(0)$ has a closed orbit $\GL(\V) \cdot \X \subset \mu_\CC^{-1}(0)$ in the quotient \eqref{eq definition of Nn_1 as a projective GIT quotient} if and only if it is stable. Furthermore, every $\X \in \mu_\CC^{-1}(0)^{st}$ has trivial stabilizer. 

Similarly, $\X \in \mu_\CC^{-1}(0)$ has a closed orbit $\GL(\V) \cdot \X \subset \mu_\CC^{-1}(0)$ in the quotient \eqref{eq definition of Nn_-1 as a projective GIT quotient} if and only if it is costable. Every $\X \in \mu_\CC^{-1}(0)^{st}$ has trivial stabilizer. 
\epr 

One has a description of $\Nn^\Qq_1$ and $\Nn^\Qq_{-1}$ in terms of simple quotients
\beq
\label{eq Nn_1 as a quotient of the stable ADHM locus}
\Nn^\Qq_1 \cong \quotient{\mu_\CC^{-1}(0)^{st}}{\GL(\V)}
\eeq
and 
\beq
\label{eq Nn_-1 as a quotient of the stable ADHM locus}
\Nn^\Qq_{-1} \cong \quotient{\mu_\CC^{-1}(0)^{ct}}{\GL(\V).}
\eeq

Recall that $\mu_\CC^{-1}(0)$ with the restriction of the symplectic form $\omega_3$ is a symplectic manifold where the action of $\U(\V)$ gives the moment map $\mu_3$. Using Kempf-Ness theorem \cite{kempf&ness} one can express the GIT quotients as symplectic quotients, that is
\beq
\label{eq definition of Nn_0 as a hyperkahler quotient}
\Nn^\Qq_0 \cong \quotient{\mu_\CC^{-1}(0) \cap \mu_3^{-1}(0)}{\U(\V)} = \quotient{\mu_1^{-1}(0) \cap \mu_2^{-1}(0) \cap \mu_3^{-1}(0)}{\U(\V)};
\eeq
taking $r \in i\RR^{>0}$, one has
\beq
\label{eq definition of Nn_1 as a hyperkahler quotient}
\Nn^\Qq_1 \cong  \quotient{\mu_\CC^{-1}(0) \cap \mu_3^{-1}(-r)}{\U(\V)} = \quotient{\mu_1^{-1}(0) \cap \mu_2^{-1}(0) \cap \mu_3^{-1}(-r)}{\U(\V)}
\eeq
and
\beq
\label{eq definition of Nn_-1 as a hyperkahler quotient}
\Nn^\Qq_{-1} \cong  \quotient{\mu_\CC^{-1}(0) \cap \mu_3^{-1}(r)}{\U(\V)} = \quotient{\mu_1^{-1}(0) \cap \mu_2^{-1}(0) \cap \mu_3^{-1}(r)}{\U(\V)}.
\eeq
Note that $\Nn^\Qq_0$, $\Nn^\Qq_1$ and $\Nn^\Qq_{-1}$ are hyperk\"ahler manifolds, defined as hyperk\"ahler quotients of $\MM^\Qq_{\V,\W}$.

Finally, we define $\Nn^\Qq_{reg}$ as the open subset of $\Nn^\Qq_1$ consisting of regular points. It is also a hyperk\"ahler manifold and recalling \eqref{eq Nn_1 as a quotient of the stable ADHM locus} and \eqref{eq Nn_-1 as a quotient of the stable ADHM locus}, one has
\beq
\label{eq Nn_reg as a quotient of the regular ADHM locus}
\Nn^\Qq_{reg} \cong \quotient{\mu_\CC^{-1}(0)^{reg}}{\GL(\V)}.
\eeq
Note that \eqref{eq definition of Nn_0 as a affine GIT quotient}, \eqref{eq definition of Nn_-1 as a projective GIT quotient}, and the fact that $\mu_\CC^{-1}(0)^{reg}$ is preserved under the action of $\GL(\V)$ imply that $\Nn^\Qq_{reg}$ can naturally be understood as an open subset of $\Nn^\Qq_0$ and $\Nn^\Qq_{-1}$.


\subsection{Branes from involutions on a hyperk\"ahler quotient}
\label{sc branes in a hyperkahler quotient}

In this section, we decribe how to construct branes on a hyperk\"ahler quotient using analytic involutions on the covering hyperk\"ahler.

To be precise, let $(Y,\Gamma_1,\Gamma_2,\Gamma_3,\eta)$ be a (finite dimensional) hyperk\"ahler manifold, let $a : Y \to Y$ be an analytic involution,  and denote by $Y^a$ its subvariety of fixed points. By abuse of notation, denote by $a : T_xY \to T_xY$ the induced involution in the tangent space of a fixed point $x\in Y^a$, and recall that the tangent space of $Y^a$ at $x$ is the invariant subspace $(T_x Y)^a$. Suppose that 
\beq \label{eq eta invariant under a}
\eta(a( \qua \cdot \qua), a( \qua \cdot \qua)) = \eta(\qua \cdot \qua, \qua \cdot \qua),
\eeq 
that is, $a$ is an isometry, and that for each $k = 1,2,3$, and  
\beq
\label{eq Gamma delta-commutes with a}
\Gamma_k a (\qua \cdot \qua) = \delta_k a \Gamma_k (\qua \cdot \qua) \qquad \text{with } \delta_k = \pm 1. 
\eeq
If $\delta_k = 1$, that is, if $a$ commutes with the complex structure $\Gamma_k$, the involution is holomorphic with respect to $\Gamma_k$. On the other hand, if $\delta_k = -1$, i.e. $a$ anticommutes with the complex structure  $\Gamma_k$ and we say that $a$ is antiholomorphic with respect to $\Gamma_k$. In the last case, one has the following result.

\blm[\cite{baraglia&schaposnik_1}, Lemma 9] \label{lm dim Y^a is one half}
Let $Y$ be a complex manifold of (complex) dimension $n$ and let $a : Y \to Y$ be an antiholomorphic involution. If the fixed point locus $Y^a$ is not empty, then it is a smooth analytic subvariety of real dimension $n$.
\elm

\brm \label{rm Y^a is isotropic}
Let us further consider $(Y,\Gamma,\eta)$ to be a K\"ahler manifold. Take an antiholomorphic involution $a : Y \to Y$ satisfying \eqref{eq eta invariant under a}. Then the fixed point locus is isotropic, that is, for every two $u,v \in T_x(Y^a) = (T_x Y)^a$, one has $\omega(u,v)= 0$, since
\begin{align*}
\omega(u,v) & = \eta(u,\Gamma v) = \eta(au,a\Gamma v)
\\
& = \eta(au, -\Gamma av) = -\eta(u,\Gamma v) 
\\
& = -\omega(u,v). 
\end{align*}
\erm 

If the involution $a$ satisfies \eqref{eq eta invariant under a} and \eqref{eq Gamma delta-commutes with a}, one has that: 
\begin{itemize}
\item  if $\delta_k = 1$ then $Y^a$ is a $B$-brane (complex subvariety) with respect to $\Gamma_k$ since the involution $a$ is holomorphic with respect to $\Gamma_k$, and  
\item  if $\delta_k = -1$ then $Y^a$ is an $A$-brane (lagrangian subvariety) with respect to $\Gamma_k$, since in that case $Y^a$ is isotropic by Remark \ref{rm Y^a is isotropic} and has maximal dimension by Lemma \ref{lm dim Y^a is one half}. 
\end{itemize}

Now let $G$ be a compact group acting on $Y$ isometrically with respect to $\eta$. One has moment maps
$\mu_1$, $\mu_2$ and $\mu_3$, each associated to a symplectic form $\omega_k$, with $\mu_k : Y \to \lie{g}^*$. For $\zeta = (\zeta_1, \zeta_2, \zeta_3)$ with $\zeta_i \in \zZ^*$ (where $\zZ$ denotes the centre of $\gG$), the associated hyperk\"ahler quotient is defined as
\[
\ol{Y}_\zeta := \quotient{\mu_1^{-1}(\zeta_1) \cap \mu_2^{-1}(\zeta_2) \cap \mu_3^{-1}(\zeta_3)}{G}.
\]
The importance of this construction relies in the fact that $\ol{Y}_\zeta$ is also a hyperk\"ahler manifold, with metric $\ol{\eta}$, complex structures $\ol{\Gamma}_1$, $\ol{\Gamma}_2$ and $\ol{\Gamma}_3$, and symplectic forms $\ol{\omega}_1$, $\ol{\omega}_2$ and $\ol{\omega}_3$. Furthermore, for the natural inclusion and projection maps of the following diagram
\[
\xymatrix{
\mu_1^{-1}(\zeta_1) \cap \mu_2^{-1}(\zeta_2) \cap \mu_3^{-1}(\zeta_3) \qua \ar@{^{(}->}[rr]^{\qquad \qquad \iota} \ar[d]^{\pi} & & Y 
\\
\ol{Y}_\zeta. & &
}
\]
we have that $\pi^* \ol{\eta}$, $\pi^*\ol{\Gamma}_k$ and $\pi^*\ol{\omega}_k$ are equal to $\iota^*\eta$, $\iota^*\Gamma_k$ and $\iota^*\omega_k$.

The involution $a : Y \to Y$ descends to an involution $\ol{a} : \ol{Y}_\zeta \to \ol{Y}_\zeta$ in the hyperk\"ahler quotient, if and only if $a$ restricts to $\mu_1^{-1}(\zeta_1) \cap \mu_2^{-1}(\zeta_2) \cap \mu_3^{-1}(\zeta_3)$, i.e. 
\beq 
\label{eq a restricts to mu^-1 0}
a\left (\mu_1^{-1}(\zeta_1) \cap \mu_2^{-1}(\zeta_2) \cap \mu_3^{-1}(\zeta_3) \right ) \subseteq \mu_1^{-1}(\zeta_1) \cap \mu_2^{-1}(\zeta_2) \cap \mu_3^{-1}(\zeta_3),
\eeq 
and, for every point $Y \in \mu_1^{-1}(\zeta_1) \cap \mu_2^{-1}(\zeta_2) \cap \mu_3^{-1}(\zeta_3)$, the image of the orbit of $Y$ is the orbit of the image of $Y$,
\beq \label{eq the image of the orbit is the orbit of the image}
a(G \cdot x) \subseteq G \cdot a(x).
\eeq 
Note that if $a$ satisfies \eqref{eq eta invariant under a} and \eqref{eq Gamma delta-commutes with a} for $\eta$ and the $\Gamma_k$, then $\ol{a}$ also satisfies \eqref{eq eta invariant under a} and \eqref{eq Gamma delta-commutes with a} for $\ol{\eta}$ and $\ol{\Gamma}_1$, $\ol{\Gamma}_2$ and $\ol{\Gamma}_3$. Therefore, we have established the main result of this section.

\blm\label{basic lemma}
Let $G$ be a compact Lie group acting isometrically on the hyperk\"ahler manifold $(Y,\Gamma_1,\Gamma_2,\Gamma_3,\eta)$, giving the hyperk\"ahler quotient $(\ol{Y},\ol{\Gamma}_1,\ol{\Gamma}_2,\ol{\Gamma}_3,\ol{\eta})$, as above. Let $a : Y \to Y$ be an involution satisfying conditions \eqref{eq eta invariant under a}, \eqref{eq Gamma delta-commutes with a}, \eqref{eq a restricts to mu^-1 0} and \eqref{eq the image of the orbit is the orbit of the image}. Then $a$ induces an involution $\ol{a} : \ol{Y}_\zeta \to \ol{Y}_\zeta$ in the quotient, and the subvariety of fixed points $(\ol{Y}_\zeta)^{\ol{a}}$ is a brane inside $\ol{Y}_\zeta$, whose type is given by the values of $\delta_1$, $\delta_2$ and $\delta_3$.
\elm


\section{Branes on Nakajima quiver varieties}

In this section we describe four different involutions on the Nakajima quiver varieties, and we show that the corresponding fixed point sets are branes of all possible types. 

First, we set up some notation. Given an involution $a:\MM^\Qq_{\V,\W}\to\MM^\Qq_{\V,\W}$, and elements $\g \in \U(\V)$ and $\h \in \U(\W)$, we define the automorphism
\beq \label{eq definition of a_GH}
\morph{a_{(\g,\h)}}{\MM^\Qq_{\V,\W}}{\MM^\Qq_{\V,\W}}{X}{(\g,\h) \cdot a(\X).}{}
\eeq
Note that $a_{(\g,\h)}$ is not always involution, but only if $\g \in \U(\V)$ and $\h \in \U(\W)$ satisfy certain conditions which will be given precisely in each case. 

Since the action of $\U(\V) \times \U(\W)$ commutes with the three complex structures $\Gamma_k$, and preserves the metric $\eta$ and the symplectic forms $\omega_k$, it is not difficult to check that if $a$ satisfies the conditions of Lemma \ref{basic lemma}, then so does $a_{(\g,\h)}$ when it is an involution.


\subsection{The transposition involution}

The first example is given by the involution 
\[
\morph{b}{\MM^\Qq_{\V,\W}}{\MM^\Qq_{\V,\W}}{(\A, \B, \I, \J)}{(\A^{t}, \B^{t}, \J^\vee, -\I^\vee).}{}
\]

\blm
\label{lm properties of b}
The involution $b$ is an isometry that commutes with $\Gamma_1$, $\Gamma_2$ and $\Gamma_3$. Furthermore, $b$ preserves $\mu^{-1}_\CC(0) \cap \mu_3^{-1}(0)$ and for every $\X \in \MM^\Qq_{\V,\W}$ one has that
\[
b \left ( \U(\V) \cdot \X \right ) = \U(\V) \cdot b(\X).
\]
\elm 

In other words, $b$ satisfies the conditions \eqref{eq eta invariant under a}, \eqref{eq Gamma delta-commutes with a} with $\delta_k=1$, \eqref{eq a restricts to mu^-1 0} with $\zeta_k=0$, and \eqref{eq the image of the orbit is the orbit of the image}.

\proof
Let $\X = (\A, \B, \I, \J)$ and $\X' = (\A', \B', \I', \J')$,
\begin{align*}
\eta(b(\X), b(\X')) = &  \frac{1}{2} \Tr \left ( \A^\vee\ol{\A}' + (\A')^\vee \ol{\A} + \B^\vee\ol{\B}' + (\B')^\vee \ol{\B}  \right ) 
\\
& \qquad + \frac{1}{2} \Tr \left ( \J^\vee \ol{\J}' + (\J')^\vee \ol{\J} \right )  + \frac{1}{2} \Tr \left ( \I^\vee\ol{\I}'  +  (\I')^\vee\ol{\I} \right ) 
\\
= & \frac{1}{2} \Tr \left (  (\A')^* \A +  \A^* \A' + (\B')^* \B +  \B^* \B' \right )^\vee 
\\
& \qquad + \frac{1}{2} \Tr \left ( (\J')^*\J + \J^*\J'  \right )^\vee  + \frac{1}{2} \Tr \left ( (\I')^*\I  + \I^*\I' \right )^\vee
\\
= & \eta(\X,\X').
\end{align*}

One trivially has that $\Gamma_1 b = b \Gamma_1$. Also, $\Gamma_2 b = b \Gamma_2$, since 
\[
\Gamma_2 b(\X) =  \left (  -\ol{\B} ,  \ol{\A}, \ol{\I},  - \ol{\J} \right ) = b \Gamma_2(\X)
\]
Commutativity with $\Gamma_3$ follows from the commutativity with $\Gamma_1$ and $\Gamma_2$.

Note that
\begin{align*}
\mu_\CC(b(\X)) = & [\A^\vee, \B^\vee] - \J^\vee \I^\vee
\\
= & - [\A,\B]^\vee -(\I\J)^\vee
\\
= & - \mu_\CC(\X)^\vee,
\end{align*}
and
\begin{align*}
\mu_3(b(\X)) = & \frac{i}{2} \left( [\A^\vee, (\A^*)^\vee ] + [\B^\vee, (\B^*)^\vee ] + \J^\vee (\J^*)^\vee - (\I^*)^\vee \I^\vee \right )
\\
= & \frac{i}{2} \left( - [\A, \A^* ]^\vee  - [\B, \B^*]^\vee +  (\J^*\J)^\vee - (\I\I^*)^\vee \right )
\\
= & - \mu_3(\X)^\vee.
\end{align*}

Finally, for every $\k \in \U(\V)$, one has that
\begin{align*}
b(\k \cdot \X) & = b \left ( \k \A \k^{-1}, \k \B \k^{-1}, \k\I, \J\k^{-1} \right )
\\
& = \left ( (\k^{-1})^\vee \A^\vee \k^\vee, (\k^{-1})^\vee \B^\vee \k^\vee, (\k^{-1})^\vee \J^\vee, - \I^\vee \k^\vee    \right )
\\
& = (\k^{-1})^\vee \cdot b(\X).
\end{align*}
\qed

Next, consider the automorphism $b_{(\g,\h)}: \MM^\Qq_{\V,\W}\to \MM^\Qq_{\V,\W}$, defined as in 
(\ref{eq definition of a_GH}). Note that $(b_{(\g,\h)})$ is an involution if and only if
\beq \label{rm G and H for b_GH involution}
\g = (g_i)_{i \in \Delta} = (e^{i\xi_k} g_i^\vee)_{i \in \Delta} \quad \nr{and} \quad \h = (h_i)_{i \in \Delta} = (-e^{i\xi_k} h_i^\vee)_{i \in \Delta}
\eeq
with
\[
\xi_{\tail(a)} = - \xi_{\head(a)}
\]
for every $a \in \Sigma$. 

\begin{remark} \label{lm two types of autodual instanton}
If $\Qq$ has a closed path with an odd number of steps, then, either
\begin{itemize}
 \item $\g^\vee = \g$ and $\h^\vee = -\h$, or
 \item $\g^\vee = -\g$ and $\h^\vee = \h$.
\end{itemize}
\end{remark}

One easily checks that, under the conditions posed above, $b_{(\g,\h)}$ satisfy the conclusions of Lemma \ref{lm properties of b}.


Summing all up (Lemma \ref{basic lemma}, Lemma \ref{lm properties of b} and the observations above), we obtain the following statement.

\bco
Take $(\g,\h)$ as specified by (\ref{rm G and H for b_GH involution}), and let $(\Nn^\Qq_0)^b$ be the fixed point locus of the involution $b_{(\g,\h)}$ on $\Nn^\Qq_0$. Then $(\Nn^\Qq_0)^b$ is a $(B,B,B)$-brane within $\Nn^\Qq_0$.
\eco

Note that $b_{(\g,\h)}$ does not descend to an involution neither on $\Nn^\Qq_1$ nor on $\Nn^\Qq_{-1}$ since $b$ does not preserve neither $\mu_3^{-1}(-r)$ nor $\mu_3^{-1}(r)$ for $r \in i\RR^{>0}$. However, it follows immediately from the proof of \cite[Lemma 2.7]{nakajima_3} that $\X$ stable implies $b_{\g,\h}(\X)$ costable and vice-versa, so that $b_{(\g,\h)}$ preserves regularity and descends to an involution $\ol{b}$ on $\Nn^\Qq_{reg}$.

\bco\label{inv b reg}
The fixed point locus $(\Nn^\Qq_{reg})^{\ol{b}}$ of the involution given by $b_{(\g,\h)}$ on $\Nn^\Qq_{reg}$, is a $(B,B,B)$-brane within $\Nn^\Qq_{reg}$.
\eco


\subsection{The sign involution}

Consider a function 
\beq
\label{eq description of gamma function}
\gamma : \Sigma \sqcup \Sigma_\Delta \lra \{ 1, -1\}
\eeq
and use it to define the involution
\[
\morph{c^\gamma}{\MM^\Qq_{\V,\W}}{\MM^\Qq_{\V,\W}}{(\A, \B, \I, \J)}{\left ( (\gamma(a)A_a)_{a \in \Sigma}, (\gamma(a)B_{\wt{a}})_{a \in \Sigma}, (\gamma(i)I_i)_{i \in \Delta}, (\gamma(i)J_i)_{i \in \Delta} \right ).}{}
\]
We assume that $\gamma$ is not the constant function $\gamma=1$, so that $c^\gamma\ne\id$.

\blm
\label{lm properties of c}
For any $\gamma$ as in \eqref{eq description of gamma function}, the involution $c^\gamma$ is an isometry that commutes with $\Gamma_1$, $\Gamma_2$ and $\Gamma_3$ and preserves $\mu^{-1}_\CC(0) \cap \mu_3^{-1}(r)$ for every $r \in i\RR$. Also, for every $\X \in \MM^\Qq_{\V,\W}$ one has that
\[
c^\gamma \left ( \U(\V) \cdot \X \right ) = \U(\V) \cdot c^\gamma(\X).
\]
\elm 

In other words, $c^\gamma$ satisfies the conditions \eqref{eq eta invariant under a}, \eqref{eq Gamma delta-commutes with a} with $\delta_k=1$, \eqref{eq a restricts to mu^-1 0} with $\zeta_1=\zeta_2=0$, $\zeta_3=ir\id$, and \eqref{eq the image of the orbit is the orbit of the image}.

\proof
Consider $\F \in \Hom^\Sigma(\V, \V)$ (resp. $\Hom^\Delta(\V, \W)$) and $\E \in \Hom^{\Sigma^\vee}(\V, \V)$ (resp. $\Hom^{\Delta}(\V,\W)$). Setting
\[
\F' = \left ( \gamma(a) F_a \right )_{a \in \Sigma} \quad 
\left (\nr{resp. }  \left (  \gamma(i) F_i \right )_{i \in \Delta} \right ),  
\]
and
\[
\E' = \left (  \gamma(a) E_{\wt{a}} \right )_{a \in \Sigma} 
\quad \left (\nr{resp. }  \left (  \gamma(i) E_i \right )_{i \in \Delta} \right ),  
\]
one can easily check that $\E'\F' = \E\F$. Then, it follows that 
\begin{itemize}
 \item $\eta(c^\gamma(\X), c^\gamma(\X')) = \eta(\X, \X')$,
 \item $\mu_3(c^\gamma(\X)) = \mu_3(\X)$, and
 \item $\mu_\CC(c^\gamma(\X)) = \mu_\CC(\X)$.
\end{itemize}

It is trivial that $\Gamma_k c^\gamma = c^\gamma \Gamma_k$ for $k = 1,2,3$, and also that for every $\k \in \U(\V)$, one has that
\[
c^\gamma \left ( \k \A \k^{-1}, \k \B \k^{-1}, \k\I, \J\k^{-1} \right ) = \left ( \k \A' \k^{-1}, \k \B' \k^{-1}, \k\I', \J'\k^{-1} \right ),
\]
so $c^\gamma(\k \cdot \X) = \k \cdot c^\gamma(\X)$.
\qed

Given $\g \in \U(\V)$ and $\h \in \U(\W)$, consider $c^\gamma_{(\g,\h)}$ as defined in \eqref{eq definition of a_GH}. Note that $(c^\gamma_{(\g,\h)})^2 = \id$ if and only if
\begin{equation} \label{rm G and H for c_GH involution}
\g^2 = (e^{i\xi_k})_{i \in \Delta} \quad \nr{and} \quad \h^2 = (e^{-i\xi_k})_{i \in \Delta}
\end{equation}
with
\[ 
\xi_{\tail(a)} = \xi_{\head(a)} 
\]
for every $a \in \Sigma$. If $\Qq$ has a closed path with an odd number of steps, then
\[
\g^2 = (e^{i\xi})_{i \in \Delta}  \quad \nr{and} \quad \h^2 = (e^{-i\xi})_{i \in \Delta}.
\]

As observed in the beginning of this section, Lemma \ref{lm properties of c} extends for $c^\gamma_{(\g,\h)}$.  It is therefore clear that $c^\gamma_{\g,\h}$ descends to an involution $\ol{c}$ on the hyperk\"ahler quotients $\Nn^\Qq_0$, $\Nn^\Qq_1$ and $\Nn^\Qq_{-1}$. Moreover, it is automatic to check that $c^\gamma_{(\g,\h)}$ preserves stability and costability;  therefore it restricts to an involution on $\Nn^\Qq_{reg}$.


\bco \label{cor brane c_GH}
Take $\gamma$ as in \eqref{eq description of gamma function}, and $(\g,\h)$ as in \eqref{rm G and H for c_GH involution}. For $* = 0, 1, -1$ and $reg$, let $(\Nn^\Qq_*)^{\ol{c}}$ be the fixed point locus given by $c^\gamma_{(\g,\h)}$. Then $(\Nn^\Qq_*)^{\ol{c}}$ is a $(B,B,B)$-brane within $\Nn^\Qq_*$.
\eco 


\subsection{The recombination involution}
Recall that, except for $\id$ and $-\id$, the idempotent elements of $\U(2)$ lie in the orbit of
$\begin{pmatrix} 1 & 0 \\ 0 & -1 \end{pmatrix}$ under the conjugation action on $\U(2)$, that is the set 
\beq
\label{eq non-trivial involution in U2}
T:= \left \{ \begin{pmatrix}
           t & z \\
           \ol{z} & -t
        \end{pmatrix}, \nr{ satisfying } t^2 +  z\ol{z} = 1, \nr{ with } t \in \RR \nr{ and } z \in \CC \right \}.
\eeq

Write $\Sigma = \Lambda \sqcup \Lambda^\perp$, where $\Lambda$ is the {\it set of loops} of the quiver,
\[
\Lambda := \{ a \in \Sigma \nr{ such that } \tail(a) = \head(a)  \},
\]
and $\Lambda^\perp$ is its complement, 
\[
\Lambda^\perp := \{ a \in \Sigma \nr{ such that } \tail(a) \neq \head(a)  \}.
\]

Next, consider a function
\beq \label{eq description of delta function}
\delta : \Lambda \sqcup \Lambda^\perp \sqcup \Sigma_\Delta \lra T, 
\eeq
such that 
$$ \delta(a) = \begin{pmatrix} t_a & 0 \\ 0 & -t_a \end{pmatrix} $$
for every $a \in \Lambda^\perp \sqcup \Sigma_\Delta$. Note that
\[
\delta(a) \begin{pmatrix} A_a \\ B_{\wt{a}} \end{pmatrix} = \begin{pmatrix}
           t_a & z_a \\
           \ol{z}_a & -t_a
        \end{pmatrix} \begin{pmatrix} A_a \\ B_{\wt{a}} \end{pmatrix} = \begin{pmatrix}
           t_a A_a + z_a B_{\wt{a}} \\
           \ol{z}_a A_a - t_a B_{\wt{a}}
        \end{pmatrix}
\]
and
\[
\delta(i) \begin{pmatrix} I_i \\ J_i \end{pmatrix} = \begin{pmatrix}
           t_i & 0 \\
           0 & -t_i
        \end{pmatrix} \begin{pmatrix} I_i \\ J_i \end{pmatrix} = \begin{pmatrix}
           t_i I_i \\
           - t_i J_i
        \end{pmatrix}
\]
Accordingly, we denote
\[
\delta_1(\A, \B) = (t_a A_a + z_a B_{\wt{a}})_{a \in \Sigma}, \quad \nr{and} \quad \delta_2(\A, \B) = (\ol{z}_a A_a - t_a B_{\wt{a}})_{a \in \Sigma}
\]
and 
\[
\delta(\I) = (t_i I_i)_{i \in \Delta}, \quad \nr{and} \quad \delta(\J) = (- t_i J_i)_{i \in \Delta}.
\]

Given $\delta$ as in \eqref{eq description of delta function}, define the following involution on $\MM^\Qq_{\V,\W}$,
\[
\morph{d^\delta}{\MM^\Qq_{\V,\W}}{\MM^\Qq_{\V,\W}}{(\A, \B, \I, \J)}{\left ( (\delta_1(\A, \B), \delta_2(\A, \B), \delta(\I), \delta(\J) \right ).}{}
\]

\blm
\label{lm properties of d}
The involution $d^\delta$ commutes with $\Gamma_1$ and anticommutes with $\Gamma_2$ and $\Gamma_3$, it is an isometry that preserves $\mu^{-1}_\CC(0) \cap \mu_3^{-1}(r)$ for every $r \in i\RR$ and, for every $\X \in \MM^\Qq_{\V,\W}$, one has that
\[
d^\delta \left ( \U(\V) \cdot \X \right ) = \U(\V) \cdot d^\delta(\X).
\] 
\elm 

In other words, $d^\delta$ satisfies the conditions \eqref{eq eta invariant under a}, \eqref{eq Gamma delta-commutes with a} with $\delta_1=1$ and $\delta_2=\delta_3=-1$, \eqref{eq a restricts to mu^-1 0} with $\zeta_1=\zeta_2=0$, $\zeta_3=ir\id$, and \eqref{eq the image of the orbit is the orbit of the image}.

\proof
Consider $\F, \F' \in \Hom^\Sigma(\V, \V)$ (resp. $\Hom^\Delta(\V, \W)$) and $\E,\E' \in \Hom^{\Sigma^\vee}(\V, \V)$ (resp. $\Hom^{\Delta}(\V,\W)$). Thanks to the fact that $\delta(a) \in \U(2)$, one can prove that 
\begin{align*}
\delta_1(\E, \F)^* \delta_1(\E',\F') + \delta_2(\E', \F') \delta_2(\E,\F)^* & =  \left ( \sum_{\tail(a)=i} \begin{pmatrix} E_a^* & F_a^*  \end{pmatrix} \delta(a)^* \delta(a) \begin{pmatrix} E'_a \\ F'_a  \end{pmatrix} \right )_{i \in \Delta}
\\
& =  \left ( \sum_{\tail(a)=i} \begin{pmatrix} E_a^* & F_a^*  \end{pmatrix} \begin{pmatrix} E'_a \\ F'_a  \end{pmatrix} \right )_{i \in \Delta}
\\
& = \E^* \E' + \F^* \F'.
\end{align*}
Then, we automatically have 
\begin{itemize}
 \item $\eta(d^\delta(\X), d^\delta(\X')) = \eta(\X, \X')$, and
 \item $\mu_3(d^\delta(\X)) = \mu_3(\X)$.
\end{itemize}

Note that $\Gamma_2 d^\delta = -d^\delta \Gamma_2$, given by 
\[
\delta_1(-\B^*, \A^*) = (-t_a B^*_{\wt{a}} + z_a A^*_a)_{a \in \Sigma} =  (-t_a B_{\wt{a}} + \ol{z}_a A_a)^*_{a \in \Sigma} = \delta_2(\A, \B)^* 
\]
and
\[
\delta_2(-\B^*, \A^*) = (-\ol{z}_a B^*_{\wt{a}} - t_a A^*_a)_{a \in \Sigma} =  -(t_a A_a + z_a B_{\wt{a}})^*_{a \in \Sigma} = -\delta_1(\A, \B)^* 
\]
Since $\Gamma_1 d^\delta = d^\delta \Gamma_1$, then $d^\delta$ anticommutes with $\Gamma_3 = \Gamma_1 \Gamma_2$.

We can check that
\begin{align*}
[t_a A_a + z_a B_{\wt{a}}, \ol{z}_a A_a - t_a B_{\wt{a}}] = & \qua t_a \ol{z_a} A^2_a - t_a^2 A_a B_{\wt{a}} +  |z_a|^2 B_{\wt{a}} A_a - z_a t_a B_{\wt{a}}^2 -
\\
& - \ol{z_a} t_a A_a^2 - |z_a|^2 A_a B_{\wt{a}}  + t_a^2 B_{\wt{a}} A_a + t_a z_a B_{\wt{a}}^2
\\
= &  - (t_a^2 + |z_a|^2) A_a B_{\wt{a}}  + (t_a^2 + |z_a|^2) B_{\wt{a}} A_a 
\\
= & - [A_a, B_{\wt{a}}],
\end{align*}
and therefore we know that
\begin{align*}
\mu_\CC(d^\delta(\X)) = & [ \delta_1(\A, \B), \delta_2(\A,\B) ] + \delta(\I) \delta(\J) 
\\
= & \left ( \sum_{\tail(a) = i} \left [ t_a A_a + z_a B_{\wt{a}}, \ol{z}_a A_a - t_a B_{\wt{a}} \right ] \right )_{i \in \Delta} + (-t_i^2 I_i J_i )_{i \in \Delta}
\\
= & \left ( \sum_{\tail(a) = i} - \left [ A_a,B_{\wt{a}} \right ] \right )_{i \in \Delta} - (I_i J_i )_{i \in \Delta}
\\
= & - \mu_\CC(\X).
\end{align*}

Finally, for every $\k \in \U(\V)$, one has that
\begin{align*}
d^\delta(\k \cdot \X) & = d^\delta \left ( \k \A \k^{-1}, \k \B \k^{-1}, \k\I, \J\k^{-1} \right )
\\
& =  \left ( \delta_1(\k \A \k^{-1}, \k \B \k^{-1}), \delta_2(\k \A \k^{-1}, \k \B \k^{-1}), \delta(\k\I), \delta(\J\k^{-1}) \right )
\\
& = \left ( \k \delta_1(\A, \B) \k^{-1}, \k \delta_2(\A, \B) \k^{-1}, \k \delta(\I), \delta(\J) \k^{-1} \right ), 
\\
& = \k \cdot d^\delta(\X).
\end{align*}
\qed

Given $\g \in \U(\V)$ and $\h \in \U(\W)$, consider $d^\delta_{(\g,\h)}$ as defined in \eqref{eq definition of a_GH}. Note that $(d^\delta_{(\g,\h)})^2 = \id$ if and only if
\beq \label{rm G and H for d_GH involution}
\g^2 = (e^{i\xi_k})_{i \in \Delta} \quad \nr{and} \quad \h^2 = (e^{-i\xi_k})_{i \in \Delta}
\eeq
with
\[ 
\xi_{\tail(a)} = \xi_{\head(a)} 
\]
for every $a \in \Sigma$. In addition, if $\Qq$ has a closed path with an odd number of steps, then
\[
\g^2 = (e^{i\xi})_{i \in \Delta}  \quad \nr{and} \quad \h^2 = (e^{-i\xi})_{i \in \Delta}.
\]

As observed in the beginning of this section, Lemma \ref{lm properties of d} extends for
$d^\delta_{(\g,\h)}$. It is therefore clear that $d^\delta_{\g,\h}$ descends to an involution $\ol{d}$ on the hyperk\"ahler quotients $\Nn^\Qq_0$, $\Nn^\Qq_1$ and $\Nn^\Qq_{-1}$. Moreover, it is straightforward to check that $d^\delta_{(\g,\h)}$ preserves stability and costability; therefore it restricts to an involution on $\Nn^\Qq_{reg}$. Note also that $d^\delta_{(\g,\h)}$ is an antiholomorphic holomorphic involution after a hyperk\"ahler rotation on $\Nn^\Qq_*$. Then, by Lemma \ref{lm dim Y^a is one half}, $(\Nn^\Qq_*)^{\ol{d}}$ is smooth and half-dimensional.

\bco \label{cor brane d_GH}
Take $\delta$ as in \eqref{eq description of delta function} and $(\g,\h)$ as in equation \eqref{rm G and H for d_GH involution}. For $* = 0, 1, -1$ and $reg$, let $(\Nn^\Qq_*)^{\ol{d}}$ be the fixed point locus given by $d^\delta_{(\g,\h)}$. Then $(\Nn^\Qq_*)^{\ol{d}}$ is a $(B,A,A)$-brane within $\Nn^\Qq_*$. If non-empty, $(\Nn^\Qq_*)^{\ol{d}}$ is a smooth submanifold of dimension $\dim_\CC(\Nn^\Qq_*)^{\ol{d}} = \frac{1}{2} \dim_\CC\Nn^\Qq_*$.
\eco

\brm
Consider the $\CC^*$-action on $\Nn_0^\Qq$, given by 
\[
\lambda \cdot (\A, \B, \I, \J) = (\A, \lambda \B, \I, \lambda \J). 
\]
It is shown in \cite[Section 5.4]{ginzburg}, that the subvariety of points fixed under this action $(\Nn_0^\Qq)^{\CC^*}$ is a complex lagrangian subvariety for $\Gamma_1$. Note that this agrees with our description, since $(\Nn_0^\Qq)^{\CC^*} \subset (\Nn_0^\Qq)^{d^\delta}$ when 
\[
\delta(a) = \bpm 1 & 0 \\ 0  & -1 \epm º
\]
for all $a \in \Lambda \sqcup \Lambda^\perp \sqcup \Sigma_\Delta$.
\erm 


\subsection{The conjugation involution}

Finally, we consider the following involution in $\MM^\Qq_{\V,\W}$
\[
\morph{e}{\MM^\Qq_{\V,\W}}{\MM^\Qq_{\V,\W}}{(\A, \B, \I, \J)}{(\ol{\A}, \ol{\B}, \ol{\I},\ol{\J}).}{}
\]

As before, we establish the basic properties of $e$.

\blm
\label{lm properties of e}
The involution $e$ is an isometry that commutes with $\Gamma_2$ and anticommutes with $\Gamma_1$ and $\Gamma_3$. Furthermore, $e$ preserves $\mu^{-1}_\CC(0) \cap \mu_3^{-1}(0)$ and for every $\X \in \MM^\Qq_{\V,\W}$ one has that
\[
e \left ( \U(\V) \cdot \X \right ) = \U(\V) \cdot e(\X).
\]
\elm 

In other words, $e$ satisfies the conditions \eqref{eq eta invariant under a}, \eqref{eq Gamma delta-commutes with a} with $\delta_2=1$ and $\delta_1=\delta_3=-1$, \eqref{eq a restricts to mu^-1 0} with $\zeta_1=\zeta_2=\zeta_3=0$, and \eqref{eq the image of the orbit is the orbit of the image}.

\proof
Note that $e \Gamma_1(\X) = e (i\X) = -i e(\X) = - e\Gamma_1(\X)$. On the other hand one can easily check that $\Gamma_2 e = e \Gamma_2$. The involution $e$ preserves the metric $\eta$, since for every two $\X = (\A, \B, \I, \J)$ and $\X' = (\A', \B', \I', \J')$, one has
\begin{align*}
\eta(e(\X), e(\X')) = &  \frac{1}{2} \Tr \left (  \ol{\A} (\A')^\vee +  \ol{\A'} \A^\vee +  \ol{\B} (\B')^\vee +  \ol{\B'} \B^\vee \right ) 
\\
& \qquad + \frac{1}{2} \Tr \left ( \ol{\I}(\I')^\vee  + \I' \I^\vee  +  \ol{\J} (\J')^\vee + \ol{\J'}\J^\vee \right ) + 
\\
= & \frac{1}{2} \Tr \left (  (\A' \A^*)^\vee +  ((\A')^* \A)^\vee + (\B' \B^*)^\vee +  ((\B')^* \B)^\vee \right ) 
\\
& \qquad + \frac{1}{2} \Tr \left ( (\I'\I^*)^\vee  + (\I(\I')^*)^\vee  +  (\J'\J^*)^\vee + (\J(\J')^*)^\vee \right )
\\
= & \frac{1}{2} \Tr \left (  \A' \A^* +  (\A')^* \A + \B' \B^* +  (\B')^* \B \right ) 
\\
& \qquad + \frac{1}{2} \Tr \left ( \I'\I^*  + \I(\I')^*  +  \J'\J^* + \J(\J')^* \right )
\\
= & \eta(\X,\X').
\end{align*}

It is also possible to show that the involution $e$ preserves $\mu^{-1}_1(0)$ and $\mu^{-1}_\CC(0)$. Observe that
\[
\mu_\CC(e(\X)) =  [\ol{\A}, \ol{\B}] - \ol{\I}\ol{\J} = \ol{\mu_\CC(\X)}
\]
and
\[
\mu_3(e(\X)) = \frac{i}{2} \left( [\ol{\A}, \ol{\A}^*] + [\ol{\B}, \ol{\B}^*] + \ol{\I}\ol{\I}^* - \ol{\J}^*\ol{\J} \right ) = \ol{\mu_3(\X)}.
\]

The last statement follows from the observation that $\ol{\k} \in \U(\V)$ provided that $\k \in \U(\V)$, since
\begin{align*}
e(\k \cdot \X) & = e \left ( \k \A \k^{-1}, \k \B \k^{-1}, \k\I, \J\k^{-1} \right )
\\
& = \left ( \ol{\k} \A \ol{\k}^{-1}, \ol{\k} \B \ol{\k}^{-1}, \ol{\k}\I, \J\ol{\k}^{-1}   \right )
\\
& = \ol{\k} \cdot e(\X).
\end{align*}
\qed

Given $\g \in \U(\V)$ and $\h \in \U(\W)$, consider $e_{(\g,\h)}$ as defined in \eqref{eq definition of a_GH}. Note that $(e_{(\g,\h)})^2 = \id$ if and only if
\beq \label{rm G and H for e_GH involution}
\g = (g_i)_{i \in \Delta} = (e^{i\xi_k} g_i^\vee)_{i \in \Delta} \quad
\nr{and} \quad \h = (h_i)_{i \in \Delta} = (e^{i\xi_k} h_i^\vee)_{i \in \Delta}
\eeq
with
\[ 
\xi_{\tail(a)} = \xi_{\head(a)} 
\]
for every $a \in \Sigma$. In addition, if $\Qq$ has a closed path with an odd number of steps, then
\[
\g^2 = (e^{i\xi})_{i \in \Delta}  \quad \nr{and} \quad \h^2 = (e^{-i\xi})_{i \in \Delta}.
\]

As observed in the beginning of this section, Lemma \ref{lm properties of e} extends for
$e_{(\g,\h)}$. It is therefore clear that $e_{\g,\h}$ descends to an involution $\ol{e}$ on the hyperk\"ahler quotient $\Nn^\Qq_0$. Note that $e_{(\g,\h)}$ does not preserve neither $\mu_3^{-1}(-r)$ nor $\mu_3^{-1}(r)$ for $r \in i\RR^{>0}$ and as a consequence, it does not descend to an involution neither in $\Nn^\Qq_1$ nor in $\Nn^\Qq_{-1}$. On the other hand, one can check (see the proof \cite[Lemma 2.7]{nakajima_3} for instance) that $e_{(\g,\h)}$ preserves regularity and therefore descends to an involution on $\Nn^\Qq_{reg}$ that we still denote by $\ol{e}$. Recall, also, that $e_{(\g,\h)}$ is an antiholomorphic involution. Then, by Lemma \ref{lm dim Y^a is one half}, $(\Nn^\Qq_*)^{\ol{e}}$ is smooth and of dimension one half. 

\bco \label{cor brane e_GH}
Take $(\g,\h)$ as in equation \eqref{rm G and H for e_GH involution}. For $* = 0$ and $reg$, let
the fixed point locus $(\Nn^\Qq_*)^{\ol{e}}$ of the involution given by $e_{(\g,\h)}$ on $\Nn^\Qq_*$. Then  $(\Nn^\Qq_*)^{\ol{e}}$ is a $(A,B,A)$-brane within $\Nn^\Qq_*$. If non-empty, $(\Nn^\Qq_*)^{\ol{e}}$ is a smooth submanifold of dimension $\dim_\CC((\Nn^\Qq_*)^{\ol{e}}) = \frac{1}{2} \dim_\CC(\Nn^\Qq_*)$. 
\eco


\subsection{Combining the involutions} \label{rm any combination gives a brane}

The involutions $b$, $c^\gamma$ and $e$ commute between them. The involution $d^\delta$ always commutes with $c^\gamma$ but only commutes with $b$ and $e$ for choices of $\delta$ that satisfy $z_a \in \RR$ for all $a \in \Lambda^\perp \sqcup \Sigma_\Delta$. Whenever they commute, any possible combination of these involutions gives another involution that preserves regularity and automatically satisfies the conditions \eqref{eq eta invariant under a}, \eqref{eq Gamma delta-commutes with a}, \eqref{eq a restricts to mu^-1 0} and \eqref{eq the image of the orbit is the orbit of the image} associated to the hyperk\"ahler quotient $\Nn^\Qq_0$ described in \eqref{eq definition of Nn_0 as a hyperkahler quotient}. Therefore, any composition of these involutions gives an involution on the Nakajima quiver varieties $\Nn^\Qq_0$ and $\Nn^\Qq_{reg}$ whose fixed point locus is a brane. 

In particular, $(A,A,B)$-branes can be obtained by considering fixed point loci of involutions of the form $(ed^\delta)_{(\g,\h)}$ and $(ebd^\delta)_{(\g,\h)}$, with the choice of $\delta$ that we specified above. More generally, the table below summarizes all the possibilities we have considered.

\begin{center} \begin{table}[H]
\begin{tabular}{| c | c | c |} \hline
 & Equation for $\g$ and $\h$   & Type of brane
\\ \hline
$b_{(\g,\h)}$ & $\g = \pm \g^\vee$ and $\h = \mp \h^\vee$  & $(B,B,B)$   \\
\hline 
$c^\gamma_{(\g,\h)}$ & $\g^2 = (e^{i\xi})_\Delta$ and $\h^2 = (e^{-i\xi})_\Delta$  & $(B,B,B)$   \\
\hline 
$d^\delta_{(\g,\h)}$ & $\g^2 = (e^{i\xi})_\Delta$ and $\h^2 = (e^{-i\xi})_\Delta$  & $(B,A,A)$   \\
\hline 
$e_{(\g,\h)}$ & $\g = \pm \g^\vee$ and $\h = \pm \h^\vee$  & $(A,B,A)$   \\
\hline 
$(eb)_{(\g,\h)}$ & $\g^2 = (e^{i\xi})_\Delta$ and $\h^2 = (e^{-i\xi})_\Delta$  & $(A,B,A)$ \\
\hline 
$(ec^\gamma)_{(\g,\h)}$ & $\g = \pm \g^\vee$ and $\h = \pm \h^\vee$ & $(A,B,A)$ \\
\hline 
$(ebc^\gamma)_{(\g,\h)}$ & $\g^2 = (e^{i\xi})_\Delta$ and $\h^2 = (e^{-i\xi})_\Delta$  & $(A,B,A)$ \\
\hline 
$(ed^\delta)_{(\g,\h)}$ & $\g = \pm \g^\vee$ and $\h = \pm \h^\vee$ & $(A,A,B)$ \\
\hline 
$(ebd^\delta)_{(\g,\h)}$  & $\g^2 = (e^{i\xi})_\Delta$ and $\h^2 = (e^{-i\xi})_\Delta$  & $(A,A,B)$ \\ 
\hline \end{tabular}
\caption{Brane types for composed involutions} \label{table1}
\end{table}\end{center}

However, it is not clear whether the fixed point set of any of the involutions we described in this section is actually non-empty. In the next section we consider a simple but already very interesting example: the Nakajima quiver variety provided by the Jordan quiver. For this case, we provide examples where these branes are non-empty.


\section{Branes in the moduli space of framed sheaves on $\PP^2$}

We denote by $\Jj$ the {\it Jordan quiver}, that is, a single vertex and a single edge-loop,
\begin{center}
\begin{tikzpicture}[every node/.style={},thick]
  \node(1) at (0,0){$\bullet$};
  \path [every node/.style={font=\sffamily\small}] 
  (1) edge [loop above] node {} (1);
\end{tikzpicture},
\end{center}
and note that $\wt{\Jj^\diamond}$ is the {\it ADHM quiver}, 
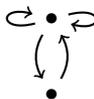
\begin{center}
\begin{tikzpicture}[every node/.style={},thick]
  \node(2) at (0,0){$\bullet$};
  \node(1) at (0,1){$\bullet$};
  \path[->]
    (1) edge [loop right] node[right] {} (1)  
        edge [loop left] node[left] {} (1)
        edge [bend right] node {} (2)
    (2) edge [bend right] node[right] {} (1);
\end{tikzpicture}.
\end{center}

Since $\Jj$ has a single vertex, one has that $\V = V$ and $\W = W$ are just vector spaces of dimensions, say, $n$ and $r$, respectively, and 
\[
\MM^\Jj_{V, W} = \Hom(V,V) \oplus \Hom(V,V) \oplus \Hom(W,V) \oplus \Hom(V,W).
\]
A point $X = (A, B, I, J)$ of $\MM^\Jj_{V, W}$ is said to satisfy the {\it ADHM equation} if 
\beq
\label{eq ADHM equation}
[A,B] + IJ = 0.
\eeq
The closed subvariety of points of $\MM^\Jj_{V,W}$ that satisfy \eqref{eq ADHM equation} is called the {\it ADHM variety} and denoted by $\AA_{V,W}$.

Note that the preimage $\mu_\CC^{-1}(0) \subset \MM_{V,W}$ is precisely the ADHM locus $\AA_{V,W}$. The definitions of stability, costability and regularity given in Section \ref{sc hyperkahler structure of Nn} apply directly in this case, and so, we have the stable and regular locus, denoted respectively by $\AA^{st}_{V,W}$ and $\AA^{reg}_{V,W}$. The action of $\U(V)$ on $\AA_{V,W}$ can be trivially extended to an action of the whole general linear group $\GL(V)$, preserving $\AA^{st}_{V,W}$ and $\AA^{reg}_{V,W}$.

On the other hand, fix a line $\ell\subset\PP^2$; recall that a \emph{framed torsion free sheaf} on $\PP^2$ is a pair $(\Ee, \Phi)$ consisting of a sheaf torsion free $\Ee$ of rank $r$, locally free on a neighborhood of $\ell$, and an isomorphism $\Phi : \Ee|_\ell \stackrel{\cong}{\lra} \Oo_\ell^{\oplus r}$, called a \emph{framing at} $\ell$.

Let $\Mm_{\PP^2}(r,n)$ denote the moduli space of framed torsion free sheaves on $\PP^2$ of rank $r$ and second Chern class, or \emph{charge}, $n$. Let also $\Mm_{\PP^2}^{\rm lf}(r,n)$ be the open subset consisting of locally free framed sheaves; recall that $\Mm_{\PP^2}^{\rm lf}(r,n)$ may also be regarded as the moduli space of framed $\SU(r)$-instantons on $\RR^4$ of charge $n$ \cite{D}.

The following is proved by Nakajima in \cite[Section 2]{nakajima_3}.

\btm\label{corresp}
There is a $1$-$1$-correspondence between 
\begin{enumerate}
\item $\GL(V)$-orbits of $X \in \AA^{st}_{V,W}$, and
\item isomorphism classes of framed torsion free sheaves with rank $r = \dim_\CC(W)$ and charge $n = \dim_\CC(V)$. 
\end{enumerate}
\etm

Recalling \eqref{eq Nn_1 as a quotient of the stable ADHM locus}, there is an isomorphism
\[
\Nn^\Jj_1 \stackrel{iso}{\cong} \Mm_{\PP^2}(r,n),
\]
that restricts to the regular locus
\[
\Nn^\Jj_{reg} \stackrel{iso}{\cong} \Mm^{\rm lf}_{\PP^2}(r,n).
\]
After these identifications, one has the following expression for the dimension,
\[
\dim_\CC(\Mm_{\PP^2}(r,n)) = \dim_\CC(\Nn^\Jj_1) = 2 rn.
\]

The goal of this section is to interpret the branes described in the previous section under the correspondence of Theorem \ref{corresp}, and provide some explicit examples showing that certain fixed point sets are nonempty. 

However, we must first recall the construction of framed torsion free sheaves on $\PP^2$ out of representation of the ADHM quiver.


\subsection{The monad construction}

Choose homogeneous coordinates $[x_0 : x_1 : x_2]$ on $\PP^2$ and set $\ell$ to be the line \{$x_0 = 0\}$. For every $X = (A, B, I, J) \in \AA_{V,W}$ define the sheaf complex
\beq \label{eq monad}
\Oo_{\PP^2}(-1) \otimes V \stackrel{\alpha^X}{\lra} 
\Oo_{\PP^2} \otimes (V \oplus V \oplus W) \stackrel{\beta^X}{\lra} 
\Oo_{\PP^2}(1) \otimes V,
\eeq
by setting
\beq
\label{eq alpha^X}
\alpha^X = \begin{pmatrix}
            x_0 A - x_1
            \\
            x_0 B - x_2
            \\
            x_0 J
           \end{pmatrix} : \Oo_{\PP^2}(-1) \otimes V \lra \Oo_{\PP^2} \otimes \begin{pmatrix}
            V
            \\
            \oplus
            \\
            V
            \\
            \oplus
            \\
            W
           \end{pmatrix}
\eeq
and
\beq
\label{eq beta^X}
\beta^X = \begin{pmatrix}
            -x_0 B + x_2, x_0 A - x_1, x_0 I
           \end{pmatrix} : \Oo_{\PP^2} \otimes \begin{pmatrix}
            V
            \\
            \oplus
            \\
            V
            \\
            \oplus
            \\
            W
           \end{pmatrix} \lra \Oo_{\PP^2}(1) \otimes V;
\eeq
note that \item $\beta^X \circ \alpha^X = 0$ if and only if $[A,B]+IJ=0$.

It is not difficult to check that (we refer to \cite[Section 2]{nakajima_3} for details): 
\begin{enumerate}
\item $\alpha^X$ is injective for any nonzero $X$,
\item $\beta^X$ is surjective if and only if $X$ is stable, in which case
$$ \Ee^X := \quotient{\ker(\beta^X)}{\im(\alpha^X)} $$
is a torsion free sheaf, and
\item $\Ee^X$ is locally free if and only if $X$ is regular. 
\end{enumerate}

Note that $\ker(\beta^X)|_\ell = (-x_1, -x_2, W)^\vee$ and $\im(\alpha^X)|_\ell = (-x_1, -x_2, 0)^\vee$, so the monad construction gives us also a framing
\[
\Phi^X : \Ee^X|_\ell = \quotient{\ker(\beta^X)|_\ell}{\im(\alpha^X)|_\ell} \stackrel{\cong}{\lra} \Oo_\ell \otimes W.
\]
The action of $M \in \GL(V)$ gives isomorphic framed sheaves, i.e. $(\Ee^X, \Phi^X) \cong (\Ee^{M \cdot X}, \Phi^{M \cdot X})$. Every framed torsion free sheaf on $\PP^2$ can be obtained in this way.


\subsection{Framed autodual sheaves}

Take $X\in\AA_{V,W}^{reg}$, and let $Y=b_{(g,h)}(X)$. It is not difficult to see that the framed locally free sheaf $\Ee^Y$ is the cohomology of the dual monad
\beq
\label{dual monad}
\Oo_{\PP^2}(-1) \otimes V^\vee \stackrel{(\beta^X)^\vee}{\lra} 
\Oo_{\PP^2} \otimes (V^\vee \oplus V^\vee \oplus W^\vee) \stackrel{(\alpha^X)^\vee}{\lra} 
\Oo_{\PP^2}(1) \otimes V^\vee,
\eeq
It follows that if $X = b_{(g,h)}(X)$, then, for any $g\in \GL(V)$ and $h\in \GL(W)$, the associated framed locally free sheaf is \emph{autodual}, i.e. there exists an isomorphism
\beq\label{iso phi}
\varphi: (\Ee^X, \Phi^X) \to \left( (\Ee^X)^\vee, (\Phi^X)^\vee \right),
\eeq
that preserves the symplectic form $h:W\otimes\Oo_\ell \to W\otimes\Oo_\ell$ on the line at infinity. Recall from \cite[Lemma 2.3]{marcos&simone&anna}, that for every $X\in\AA_{V,W}^{reg}$ satisfying $X=b_{(g,h)}(X)$, one has 
\begin{equation} \label{eq autodual monads}
X = (A, B, I, -h^{-1} I^\vee g) \textnormal{ with } gA = A^\vee g \textnormal{ and } gB = B^\vee g.
\end{equation}

As a special case of Corollary \ref{inv b reg}, we obtain the following statement. 

\bco
The submanifold consisting of framed autodual locally free sheaves is a $(B,B,B)$-brane inside
$\Mm^{\rm lf}_{\PP^2}(r,n)$.
\eco

From Remark \ref{lm two types of autodual instanton} we know that there are only two possible types of autodual locally free sheaves: \emph{symplectic} (if the isomorphism in \eqref{iso phi} satisfies $\varphi^\vee=-\varphi$) and \emph{orthogonal} bundles (if the isomorphism in \eqref{iso phi} satisfies $\varphi^\vee=\varphi$).  

For the first case, take for instance $g=\id_V$ and let $h$ be a fixed antisymmetric element of $\GL(W)$. The corresponding fixed point set $(\Nn^\Jj_{reg})^{\ol{b}}$ coincides with the moduli space of framed symplectic locally free sheaves (cf. \cite[Proposition 3.1]{marcos&simone&anna}) preserving the fixed symplectic form $h$ on the line at infinity; in other words, $(\Nn^\Jj_{reg})^{\ol{b}}$ coincides with the moduli space of framed principal ${\rm Sp}(r)$ bundles on $\PP^2$. This is known to be nonempty whenever $r$ is even; in addition, it is smooth, connected and therefore irreducible (see for instance \cite[Section 5]{scalise}). Using \eqref{eq autodual monads} one can compute the dimension as follows: the number of degrees of freedom for $A$ and $B$ is $\frac{1}{2}n(n+1)$ each, since they are symmetric; the number of degrees of freedom for $I$ is $rn$ and $J$ is fixed by $g$, $h$ and $I$. on the other hand, the ADHM equation imposes  $\frac{1}{2}n(n-1)$ independent conditions. Finally, the action of $\GL(V)$ reduces to $\Ort(V)$, which has dimension $\frac{1}{2}n(n-1)$, since we have to restrict ourselves to those $M \in \GL(V)$ satisfying $M M^t = \id_V$. We can check that the dimension of the moduli space of framed symplectic locally free sheaves is $n(r + 2)$ whenever nonempty; note that this is even, as expected for a $(B,B,B)$-brane.

For orthogonal sheaves, one can take $g = \Omega_V$ to be the standard symplectic form on $V$ and $h$ some fixed symmetric element of $\GL(W)$. The corresponding fixed point set $(\Nn^\Jj_{reg})^{\ol{b}}$ coincides with the moduli space of framed orthogonal locally free sheaves (cf. \cite[Proposition 4.1]{marcos&simone&anna}) preserving the fixed orthogonal form $h$ on the line at infinity; in other words, $(\Nn^\Jj_{reg})^{\ol{b}}$ coincides with the moduli space of framed principal ${\rm O}(r)$ bundles on $\PP^2$. In this case, framed orthogonal locally free sheaves can only exist when $n$ is even and larger than 2, cf. \cite[Lemma 4.3 and Theorem 4.4]{marcos&simone&anna}. Using again \eqref{eq autodual monads}, we obtain that the number of degrees of freedom of $A$ and $B$ is $\frac{1}{2}n(n-1)$ each, the number of degrees of freedom of $I$ is $nr$, and $J$ is fixed by $g$, $h$ and $I$. The number of linearly independent equations is $\frac{1}{2}n(n-1)$. In this case the action of $\GL(V)$ reduces to $\Sp(V)$, of dimension $\frac{1}{2}n(n+1)$, since this is the subgroup of $\GL(M)$ preserving the standard symplectic form $\Omega_V$. Therefore, the dimension of the moduli space of framed orthogonal locally free sheaves is $n(r - 1)$ whenever non-empty; again, note that this is even since $n$ is even, as is required in the case of $(B,B,B)$-branes. 



\subsection{Framed sheaves fixed by unitary involutions}

Since $\Aut(\PP^2)=\PSL(3,\CC)$, a holomorphic involution of $\PP^2$ can be described by an idempotent element of $\PSL(3,\CC)$. We say that an involution is \emph{unitary} if it is contained in $\PSU(3) \subset \PSL(3,\CC)$.
 
\blm
\label{lm classification of involutions in PP2}
Every (non-trivial) unitary involution preserving the line $\ell = \{x_0 = 0\}$ is either
\begin{itemize}
\item[(1)] the involution $\sigma_1: \PP^2\to \PP^2$ given by 
$$ [x_0:x_1:x_2] ~~ \mapsto ~~ [-x_0:x_1:x_2] $$
\end{itemize}
or
\begin{itemize}
\item[(2)] involutions $\sigma_2: \PP^2\to \PP^2$ of the form
$$ [x_0:x_1:x_2] ~~ \mapsto ~~  [x_0:t x_1 + z x_2 : \ol{z}x_1 - t x_2] $$
where $t\in\RR$ and $z\in\CC$ satisfying $t^2 + |z|^2 = 1$
\end{itemize}
\elm 

\proof
If an involution preserves the line $\ell = \{x_0 = 0\}$, then it must be associated to a $\CC^*$-class of matrices of the form
\[
\begin{pmatrix}
 \pm 1 & 0 & 0 \\
 0 & * & * \\
 0 & * & *
\end{pmatrix}.
\]
Therefore, the proof follows from Formula \eqref{eq non-trivial involution in U2}
\qed

Let us define, with an abuse of notation, for $k=1,2$,
\[
\morph{\sigma_k}{\Mm_{\PP^2}(r,n)}{\Mm_{\PP^2}(r,n)}{(\Ee, \Phi)}
{\left ( \sigma_k^*(\Ee), \sigma_k^*(\Phi) \right ).}{}
\]
In this section we show that the fixed point sets $\Mm_{\PP^2}(r,n)^{\sigma_k}$ are a $(B,B,B)$-brane for $k=1$ and a $(B,A,A)$-brane for $k=2$.

First, we analyze the relation between the points fixed by the involution
$$ \morph{c}{\MM^\Jj_{V, W}}{\MM^\Jj_{V, W}}{(A,B,I,J)}{(-A,-B,-I,-J).}{} $$
and framed sheaves invariant under pullback by $\sigma_1$. Note that
\[ \alpha^{c(X)} = \sigma_1^*\alpha^X =
\begin{pmatrix}
x_0 (-A) - x_1 \\
x_0 (-B) - x_2 \\
x_0 (-J)
\end{pmatrix}, \]
and
\[ \beta^{c(X)} = \sigma_1^*\beta^X = 
\begin{pmatrix}
-x_0 (-B) + x_2, x_0 (-A) - x_1, x_0 (-I)
\end{pmatrix}. \]
It follows that $\Ee^{c(X)}=\sigma_1^*\Ee^X$. 

\blm \label{lm framed sheaf for c_GH}
Let $X \in \AA_{V,W}^{st}$ with $X = c_{(g,h)}(X)$. Then there exists an isomorphism of framed sheaves
\[
\varphi_1 : (\Ee^X, \Phi^X) \stackrel{\cong}{\lra} \left ( \sigma_1^*(\Ee^X), \sigma_1^* (\Phi^X) \right ).
\]
Furthermore, the restriction to the line at infinity $\ell \subset \PP^2$ is
\[
\varphi_1 |_{\ell} = h.
\]
\elm 

\proof
Note that $\sigma_1^*(\Ee^X)$ is associated to the monad
\[
\Oo_{\PP^2}(-1) \otimes V \stackrel{\sigma_1^*\alpha^X}{\lra} \Oo_{\PP^2} \otimes (V \oplus V \oplus W) \stackrel{\sigma_1^*\beta^X}{\lra} \Oo_{\PP^2}(1) \otimes V.
\]
Given $g \in \U(V)$ and $h \in \U(W)$, one can set $f:V\oplus V\oplus W \to V\oplus V\oplus W$ as follows:
\[
f = \begin{pmatrix}
     g & 0 & 0\\
     0 & g & 0 \\
     0 & 0 & h
    \end{pmatrix}.
\]
If $X = c_{(g,h)}(X)$, the following isomorphism of monads
\[
\xymatrix{
\Oo_{\PP^2}(-1)  \otimes V \ar[rr]^{\alpha^X} \ar[d]_{\id \otimes g} & & \Oo_{\PP^2} \otimes (V \oplus V \oplus W) \ar[d]^{\id \otimes f} \ar[rr]^{\beta^X} & & \Oo_{\PP^2}(1) \otimes V \ar[d]^{\id \otimes g}
\\
\Oo_{\PP^2}(-1)  \otimes V \ar[rr]^{\sigma_1^*\alpha^X}  & & \Oo_{\PP^2} \otimes (V \oplus V \oplus W) \ar[rr]^{\sigma_1^*\beta^X} & & \Oo_{\PP^2}(1) \otimes V,
}
\]
provides the desired isomorphism of framed sheaves.
\qed

Together with Corollary \ref{cor brane c_GH}, we conclude:

\bco
The fixed point locus $\Mm_{\PP^2}(r,n)^{\sigma_1}$ is a $(B,B,B)$-brane in $\Mm_{\PP^2}(r,n)$; in addition, the fixed point locus $\Mm^{\rm lf}_{\PP^2}(r,n)^{\sigma_1}$ is a $(B,B,B)$-brane inside $\Mm^{\rm lf}_{\PP^2}(r,n)$.
\eco 

\begin{remark} \label{rmk c example}
Fix $r=n=2$, and consider the involution $c^\gamma:\MM^\Jj_{V, W}\to \MM^\Jj_{V, W}$ defined by
\beq\label{c example}
c^\gamma_{(g,h)} (A,B,I,J) = (-gAg^{-1}, -gBg^{-1}, gIh, h^{-1}Jg^{-1}), 
\eeq
choosing $g \in \GL(V)$ and $h \in \GL(W)$ as
$$ g = h =
\left(
\begin{array}{cc}
0 & -1 \\
1 & 0
\end{array}
\right) $$
The fixed points of $c^\gamma_{(g,h)}$ are defined by the following conditions
\begin{description}
\item[i)] $Ag + gA = 0$;
\item[ii)] $Bg + gB = 0$;
\item[iii)] $I - gIh = 0$;
\item[iv)] $J - h^{-1}Jg^{-1} = 0$.
\end{description}
Taking the matrices
$$
A=
\left(
\begin{array}{cc}
1 & 2 \\
2 & -1
\end{array}
\right),
\:\: B=
\left(
\begin{array}{cc}
1 & 1 \\
1 & -1
\end{array}
\right),
\: \: I=
\left(
\begin{array}{cc}
1 & 0 \\
0 & 1
\end{array}
\right),
\:\:J=
\left(
\begin{array}{cc}
0 & 2 \\
-2 & 0
\end{array}
\right)
$$
we get a framed locally free sheaf; indeed, it is not difficult to check that the ADHM equation and the regularity condition are satisfied. In addition, it is fixed by the involution $c^\gamma$ in \eqref{c example}. 

Therefore, we obtain a nonempty $(B,B,B)$-brane inside $\Mm^{\rm lf}_{\PP^2}(2,2)$. Moreover, note that $2k\times 2k$ matrices whose blocks form is as above, with each entry promoted to a $k\times k$ block being a multiple of the identity matrix, provides an explicit example of a point in $\Mm^{\rm lf}_{\PP^2}(2k,2k)$ fixed by the involution \eqref{c example}, and hence a nonempty $(B,B,B)$-brane inside $\Mm^{\rm lf}_{\PP^2}(2k,2k)$.

It would be interesting to determine precisely for which values of rank and charge such $(B,B,B)$-branes are nonempty, as well as checking irreducibility, smoothness, and computing its dimension.
\end{remark}

\bigskip

The analysis of the relation between the points fixed by the involution
$$ \morph{d^\delta}{\MM^\Jj_{V, W}}{\MM^\Jj_{V, W}}{(A,B,I,J)}{(tA+zB,\overline{z}A-tB,I,J),}{} $$
and framed sheaves invariant under pullback by $\sigma_2$ is similar. Note that
\[ \sigma_2^*\alpha^X =
\begin{pmatrix}
x_0 A - (t x_1 + z x_2) \\
x_0 B - (\ol{z} x_1 - t x_2) \\
x_0 J
\end{pmatrix}, \]
and
\[ \sigma_2^*\beta^X =
\begin{pmatrix}
-x_0 B + (\ol{z} x_1 - t x_2), x_0 A - (t x_1 + z x_2), x_0 I
\end{pmatrix}. \]
Defining $f:V\oplus V\oplus W \to V\oplus V\oplus W$ in the following way
\[ f = \begin{pmatrix}
t g & z g & 0\\
\ol{z} g & -t g & 0 \\
0 & 0 & h
\end{pmatrix} \]
the diagram
\[
\xymatrix{
\Oo_{\PP^2}(-1)  \otimes V \ar[rr]^{\sigma_2^*\alpha^X} \ar[d]_{\id \otimes(-g)} & & \Oo_{\PP^2} \otimes (V \oplus V \oplus W) \ar[d]^{\id \otimes f} \ar[rr]^{\sigma_2^*\beta^X} & & \Oo_{\PP^2}(1) \otimes V \ar[d]^{\id \otimes g}
\\
\Oo_{\PP^2}(-1)  \otimes V \ar[rr]^{\alpha^{d^\delta_{(g,h)}(X)}}  & & \Oo_{\PP^2} \otimes (V \oplus V \oplus W) \ar[rr]^{\beta^{d^\delta_{(g,h)}(X)}} & & \Oo_{\PP^2}(1) \otimes V,
}
\]
is an isomorphism of monads that induces the isomorphism $\sigma_2^*\Ee^X\simeq \Ee^{d^\delta_{(g,h)}(X)}$.
Together with Corollary \ref{cor brane d_GH}, we conclude:

\bco
The fixed point locus $\Mm_{\PP^2}(r,n)^{\sigma_2}$ is a $(B,A,A)$-brane in $\Mm_{\PP^2}(r,n)$; in addition, $\Mm^{\rm lf}_{\PP^2}(r,n)^{\sigma_2}$ is a $(B,A,A)$-brane in $\Mm^{\rm lf}_{\PP^2}(r,n)$. If non-empty, $\Mm^{\rm lf}_{\PP^2}(r,n)^{\sigma_2}$ is smooth and of dimension $rn$ (whenever nonempty). 
\eco

\begin{remark} \label{rmk d example}
Fix $r=n=4$, and consider the involution $bd_{(g,h)}:\MM^\Jj_{V, W}\to \MM^\Jj_{V, W}$ defined by
\beq\label{d example}
bd_{(g,h)} (A,B,I,J) = (gA^\vee g^{-1}, - gB^\vee g^{-1}, -gJ^\vee h, -h^{-1}I^\vee g^{-1}), 
\eeq
choosing $g \in \GL(V)$ and $h \in \GL(W)$ as
$$
g = 
\left(
\begin{array}{cccc}
0 & 0 & -1 & 0 \\
0 & 0 & 0 & -1 \\
1 & 0 & 0 & 0 \\
0 & 1 & 0 & 0
\end{array}
\right)
\:\:\mbox{and}\:\:
h = 
\left(
\begin{array}{cccc}
0 & 0 & 1 & 0 \\
0 & 0 & 0 & 1 \\
-1 & 0 & 0 & 0 \\
0 & -1 & 0 & 0
\end{array}
\right)
$$
Observe that the fixed points of $bd_{(g,h)}$ are defined by the following conditions
\begin{description}
\item[i)] $Ag - gA^\vee = 0$;
\item[ii)] $Bg + gB^\vee = 0$;
\item[iii)] $I - gJ^\vee h = 0$;
\item[iv)] $J - h^{-1}I^\vee g^{-1} = 0$.
\end{description}
Notice that conditions (iii) and (iv) are equivalent with the given choices of $g$ and $h$. Now consider the matrices
$$
A =
\left(
\begin{array}{cccc}
a & 0 & 0 & 1 \\
0 & a & -1 & 0 \\
0 & 1 & a & 0 \\
-1 & 0 & 0 & a
\end{array}
\right),
\:\:\:
B = 
\left(
\begin{array}{cccc}
b_1 & b_2 & b_3 & 1 \\
0 & -b_1 & 1 & b_4 \\
b_4 & 1 & -b_1 & 0 \\
1 & b_3 & -b_2 & b_1
\end{array}
\right),\:\:\:
I = 
\left(
\begin{array}{cccc}
0 & 0 & 2 & 0 \\
0 & 0 & 0 & -2 \\
1 & 0 & 0 & 0 \\
0 & 1 & 0 & 0
\end{array}
\right).
$$
Taking $J=I$, it is easy to check that $(A,B,I,J)$ are a fixed points of $bd_{(g,h)}$. Moreover, they satisfy the ADHM equation and, since $I$ and $J$ maximal rank, also fulfill the regularity condition. It follows that $(A,B,I,J)$ as above defines a framed locally free sheaf. Checking with Table \ref{table1}, we obtain a nonempty $(B,A,A)$-brane inside $\Mm^{\rm lf}_{\PP^2}(4,4)$.

Furthermore, promoting each entry of the matrices above to a multiple of the $k\times k$ the identity matrix
provides an explicit example of a point in $\Mm^{\rm lf}_{\PP^2}(4k,4k)$ fixed by the involution \eqref{d example}, and hence a nonempty $(B,A,A)$-brane inside $\Mm^{\rm lf}_{\PP^2}(4k,4k)$, smooth and of dimension $16k^2$.
\end{remark}


\subsection{Framed real sheaves}

Atiyah \cite{atiyah} introduced the notion of real vector bundle on an algebraic variety $X$ with an antiholomorphic involution $\tau : X \to X$, as a holomorphic vector bundle $E$ equipped with an involution $\wt{\tau} : E \to E$, making the diagram 
\begin{equation} \label{eq real vector bundle}
\xymatrix{
E \ar[r]^{\wt{\tau}} \ar[d] & E \ar[d]
\\
X \ar[r]^{\tau} & X,
}
\end{equation}
commutative, and such that $\wt{\tau}$ is $\CC$-antilinear on the fibres.

Given a holomorphic vector bundle $E$, we denote by $\ol{E}$ the complex vector bundle over $X$ whose underlying $C^\infty$-bundle is the same as $E$ and it is endowed with the conjugate complex structure. Note that $\tau^*\ol{E}$ is also a holomorphic vector bundle over $X$. It is well known (see for instance \cite{baraglia&schaposnik_1, baraglia&schaposnik_2, oscar&biswas, biswas&huisman&hurtubise}) that from a real vector bundle $(E,\wt{\tau})$ on gets an isomorphism of holomorphic vector bundles $\psi : E \to \tau^*\ol{E}$ and, conversely, given any such isomorphism one can construct an involution satisfying \eqref{eq real vector bundle} $\CC$-antilinear on the fibres. Therefore, a real vector bundle is equivalent to a pair $(E, \psi : E \stackrel{\cong}{\to} \tau^*\ol{E})$. 

If we denote by $\Ee$ the locally free sheaf associated to the vector bundle $E$, we set $\tau^*\ol{E}$ to be the locally free sheaf associated to $\tau^*\ol{E}$. We say that {\it a framed real locally free sheaf} is a triple $(\Ee, \Phi, \psi)$, where $(\Ee, \Phi)$ is a framed locally free sheaf and $\psi : (\Ee, \Phi) \to (\tau^*\ol{\Ee}, \tau^*\ol{\Phi})$ is an isomorphism of framed torsion free sheaves.

Let us consider the antiholomorphic involution on $\PP^2$ given by conjugation:
\[
\morph{\tau_0}{\PP^2}{\PP^2}{[x_0 : x_1 : x_2]}{[\ol{x_0} : \ol{x_1} : \ol{x_2}].}{}
\]
Since $(\tau_0^*\ol{\Ee^X}, \tau_0^*\ol{\Phi^X})$ is stable provided $(\Ee, \Phi)$ is stable, one has the following involution on the moduli space
\[
\morph{\ol{\tau}_0}{\Mm^{\rm lf}_{\PP^2}(r,n)}{\Mm^{\rm lf}_{\PP^2}(r,n)}{(\Ee, \Phi)}{(\tau_0^*\ol{\Ee^X}, \tau_0^*\ol{\Phi^X}).}{}
\]
Consequently, we refer to the fixed point set $\Mm^{\rm lf}_{\PP^2}(r,n)^{\ol{\tau}_0}$ as the locus of framed real locally free sheaves associated to $\tau_0$.




We observe that $\ol{\tau}_0$ corresponds with that involution induced by $e_{(g,h)}$.

\blm \label{lm framed sheaf for e_GH}
Let $X \in \AA_{V,W}^{reg}$ with $X = e_{(g,h)}(X)$. Then the associated framed locally free sheaf satisfies
\[
\left ( \Ee^X, \Phi^X \right ) \cong \left ( \tau_0^*\ol{\Ee^X}, \tau_0^*\ol{\Phi^X} \right ).
\]
\elm 

\proof
If $(\Ee^X,\Phi^X)$ is given by the monad \eqref{eq monad} with maps \eqref{eq alpha^X} and \eqref{eq beta^X}, then, the framed locally free sheaf of $\overline{\Oo_{\PP^2}}$-modules $\left ( \ol{\Ee^X}, \ol{\Phi^X} \right )$ is given by
\[
\ol{\Oo_{\PP^2}}(-1) \otimes \ol{V} \stackrel{\ol{\alpha^X}}{\lra} \ol{\Oo_{\PP^2}} \otimes \left ( \ol{V} \oplus \ol{V} \oplus \ol{W} \right ) \stackrel{\ol{\beta^X}}{\lra} \ol{\Oo_{\PP^2}}(1) \otimes \ol{V},
\]
where
\[
\ol{\alpha^X} = \begin{pmatrix}
            \ol{x_0} \ol{A} - \ol{x_1}
            \\
            \ol{x_0} \ol{B} - \ol{x_2}
            \\
            \ol{x_0} \ol{J}
           \end{pmatrix}
\]
and
\[
\ol{\beta^X} = \begin{pmatrix}
            -\ol{x_0} \ol{B} + \ol{x_2}, \ol{x_0} \ol{A} - \ol{x_1}, \ol{x_0} \ol{I}
           \end{pmatrix}.
\]
Then, $\tau_0^*(\ol{\Ee^X})$ is given by the monad
\[
\Oo_{\PP^2}(-1) \otimes \ol{V} \stackrel{\tau_0^*\ol{\alpha^X}}{\lra} \Oo_{\PP^2} \otimes \left ( \ol{V} \oplus \ol{V} \oplus \ol{W} \right ) \stackrel{\tau_0^*\ol{\beta^X}}{\lra} \Oo_{\PP^2}(1) \otimes \ol{V},
\]
with
\[
\tau_0^*\ol{\alpha^X} = \begin{pmatrix}
            x_0 \ol{A} - x_1
            \\
            x_0 \ol{B} - x_2
            \\
            x_0 \ol{J}
           \end{pmatrix}
\]
and
\[
\tau_0^*\ol{\beta^X} = \begin{pmatrix}
            -x_0 \ol{B} + x_2, x_0 \ol{A} - x_1, x_0 \ol{I}
           \end{pmatrix}.
\]
\qed 

Together with Corollary \ref{cor brane e_GH}, we conclude:

\bco
The locus of framed real locally free sheaves $\Mm^{\rm lf}_{\PP^2}(r,n)^{\ol{\tau}_0}$ associated to $\tau_0$ is a $(A,B,A)$-brane in $\Mm^{\rm lf}_{\PP^2}(r,n)$, smooth and of dimension $rn$ (whenever nonempty).
\eco


\subsection{Combining involutions}

Note that $\tau_0$ commutes with $\sigma_1$, $\sigma_2$, and with ta\-king duals. Thus we can consider the antiholomorphic involutions of $\PP^2$
\[
\tau_1 := \tau_0 \sigma_1 = \sigma_1 \tau_0
\]
and
\[
\tau_2 := \tau_0 \sigma_2 = \sigma_2 \tau_0.
\]

Recall that $(\Ee^*, \Phi^*)$ denotes $(\ol{\Ee}^\vee, \ol{\Phi}^\vee)$. We obtain, thanks to Section \ref{rm any combination gives a brane}, the following description of a collection of branes inside $\Mm^{reg}_{\PP^2}(r,n)$. 

\newpage

\begin{center}\begin{table}[H]
\begin{tabular}{| c | c | c | c |}\hline
 & Equation for $g$ and $h$  & Involution on $\Mm^{reg}_{\PP^2}(r,n)$ & Type of brane
\\ \hline
$b_{(g,h)}$ & $g = \pm g^\vee$ and $h = \mp h^\vee$  & $(\Ee, \Phi) \mapsto (\Ee^\lor, \Phi^\lor)$ & $(B,B,B)$   \\
\hline 
$c_{(g,h)}$ & $g^2 = e^{i\xi}$ and $h^2 = e^{-i\xi}$  & $(\Ee, \Phi) \mapsto (\sigma_1^*{\Ee}, \sigma_1^*{\Phi})$  & $(B,B,B)$   \\
\hline 
$d^\delta_{(g,h)}$ & $g^2 = e^{i\xi}$ and $h^2 = e^{-i\xi}$ & $(\Ee, \Phi) \mapsto (\sigma_2^*{\Ee}, \sigma_2^*{\Phi})$  & $(B,A,A)$   \\
 \hline
$e_{(g,h)}$ & $g = \pm g^\vee$ and $h = \pm h^\vee$ & $(\Ee, \Phi) \mapsto (\tau_0^*\ol{\Ee}, \tau_0^*\ol{\Phi})$ & $(A,B,A)$ \\
\hline 
$(eb)_{(g,h)}$ & $g^2 = e^{i\xi}$ and $h^2 = e^{-i\xi}$ & $(\Ee, \Phi) \mapsto (\tau_0^*\Ee^*, \tau_0^*\Phi^*)$ & $(A,B,A)$ \\
\hline 
$(ec)_{(g,h)}$ & $g = \pm g^\vee$ and $h = \pm h^\vee$ & $(\Ee, \Phi) \mapsto (\tau_1^*\ol{\Ee}, \tau_1^*\ol{\Phi})$ & $(A,B,A)$ \\
\hline 
$(ebc)_{(g,h)}$ & $g^2 = e^{i\xi}$ and $h^2 = e^{-i\xi}$ & $(\Ee, \Phi) \mapsto (\tau_1^*\Ee^*, \tau_1^*\Phi^*)$ & $(A,B,A)$ \\
\hline 
$(ed)_{(g,h)}$ & $g = \pm g^\vee$ and $h = \pm h^\vee$ & $(\Ee, \Phi) \mapsto (\tau_2^*\ol{\Ee}, \tau_2^*\ol{\Phi})$ & $(A,A,B)$ \\
\hline 
$(ebd)_{(g,h)}$  & $g^2 = e^{i\xi}$ and $h^2 = e^{-i\xi}$ & $(\Ee, \Phi) \mapsto (\tau_2^*\Ee^*, \tau_2^*\Phi^*)$ & $(A,A,B)$ \\ \hline
\end{tabular}
\caption{Brane types for composed involutions} \label{table2}
\end{table}\end{center}

\

\begin{remark} 
The matrices considered in examples of Remark \ref{rmk c example} are real ones, and thus are fixed points of the composed involution $ec^\gamma_{(g,h)}$, with $c^\gamma_{(g,h)}$ defined as in equation \eqref{c example}.

Therefore, we obtain a nonempty $(A,B,A)$-brane inside $\Mm^{\rm lf}_{\PP^2}(2k,2k)$, smooth and of dimension $4k^2$.
\end{remark}

\bigskip

\begin{remark} 
If the entries $a,b_1,b_2,b_3,b_4$ of the matrices in Remark \ref{rmk d example} are chosen to be real, then they are 
fixed points of the composed involution $ebd_{(g,h)}$, with $bd_{(g,h)}$ defined as in equation \eqref{d example}.

Therefore, we obtain a nonempty $(A,A,B)$-brane inside $\Mm^{\rm lf}_{\PP^2}(4k,4k)$, smooth and of dimension $16k^2$.
\end{remark}


\end{document}